\documentclass[letterpaper,12pt]{article}
\usepackage[T1]{fontenc}
\usepackage{lmodern,amsmath,amsthm,amsfonts,amssymb,graphicx,float,microtype,thmtools,underscore,mathtools,xurl,thm-restate,fontawesome5}
\usepackage[dvipsnames,svgnames,table]{xcolor}
\usepackage[shortlabels]{enumitem}
\setlist[itemize]{topsep=0ex,itemsep=0ex,parsep=0ex}
\setlist[enumerate]{topsep=0ex,itemsep=0ex,parsep=0ex}
\usepackage[unicode=true]{hyperref}
\hypersetup{
colorlinks,
breaklinks=true,
linkcolor={blue!60!black},
citecolor={black},
urlcolor={blue!60!black},
pdftitle={Treewidth 3 in the Product Structure Theorem for Bounded Degree Graphs in Minor-Closed Classes}}
\usepackage[capitalise, compress, nameinlink, noabbrev]{cleveref}
\crefname{lem}{Lemma}{Lemmas}
\crefname{thm}{Theorem}{Theorems}
\crefname{cor}{Corollary}{Corollaries} 

\newcommand{\mks}{8\max\{2,g,k,a\}} 
\newcommand{\defn}[1]{\textcolor{Maroon}{\emph{#1}}}
\newcommand{\mathdefn}[1]{\textcolor{Maroon}{#1}}
\usepackage[longnamesfirst,numbers,sort&compress]{natbib}
\makeatletter
\def\NAT@spacechar{~}
\makeatother
\usepackage[margin=1in]{geometry}
\renewcommand{\baselinestretch}{1.1}
\allowdisplaybreaks

\DeclarePairedDelimiter{\floor}{\lfloor}{\rfloor}

\renewcommand{\epsilon}{\varepsilon}
\renewcommand{\emptyset}{\varnothing}
\renewcommand{\ge}{\geqslant}
\renewcommand{\le}{\leqslant}
\renewcommand{\geq}{\geqslant}
\renewcommand{\leq}{\leqslant}
\DeclareMathOperator{\topp}{top}

\DeclareMathOperator{\dist}{dist}
\DeclareMathOperator{\col}{col}
\DeclareMathOperator{\tw}{tw}

\DeclareMathOperator{\wcol}{wcol}
\DeclareMathOperator\WReach{WReach}
\DeclareMathOperator{\apex}{apex}
\newcommand{\torso}[2]{{#1}\langle{#2}\rangle}

\newcommand{\calG}{\mathcal{G}}
\newcommand{\calI}{\mathcal{I}}

\newcommand{\calF}{\mathcal{F}}

\newcommand{\calW}{\mathcal{W}}
\newcommand{\calP}{\mathcal{P}}
\newcommand{\calQ}{\mathcal{Q}}
\newcommand{\calR}{\mathcal{R}}

\newcommand{\calT}{\mathcal{T}}
\newcommand{\calL}{\mathcal{L}}
\newcommand{\calE}{\mathcal{E}}

\newcommand{\GG}{\mathcal{G}}

\newcommand{\NN}{\mathbb{N}}

\renewcommand{\thefootnote}{\fnsymbol{footnote}}
\theoremstyle{plain}
\newtheorem{thm}{Theorem}
\newtheorem{lem}[thm]{Lemma}
\newtheorem{cor}[thm]{Corollary}
\newtheorem{obs}[thm]{Observation}

\crefname{obs}{Observation}{Observations}
\newtheorem*{lem*}{Lemma}
\theoremstyle{definition}
\newtheorem{conj}[thm]{Conjecture}
\newtheorem*{conj*}{Conjecture}

\begin{document}
\title{\bf\boldmath\fontsize{18pt}{20pt}\selectfont Treewidth 3 in the Product Structure Theorem for Bounded Degree Graphs in Minor-Closed Classes}

\author{
Vida Dujmovi\'c\,\footnotemark[2]
\qquad Robert Hickingbotham\,\footnotemark[3]
\\[0.5ex]
Pat Morin\,\footnotemark[4]
\qquad David~R.~Wood\,\footnotemark[5]
}

\newcommand{\pat}[1]{{\color{red} [PM: #1]}}
\maketitle

\begin{abstract} 
  We show that, for all positive integers $t$ and $\Delta$, there exists an integer $c$ such that, for every $K_t$-minor-free graph $G$ with maximum degree at most $\Delta$, there exists a graph $H$ of treewidth at most $3$ and a path $P$ such that $G$ is isomorphic to a subgraph of the strong product $H\boxtimes P\boxtimes K_c$.   This improves on a previous result of \citet{DEMWW22} by bounding the treewidth of $H$ independently of $t$ and $\Delta$. As an immediate consequence we obtain a bound of $O(r^4)$ on the $r$-weak colouring number of bounded-degree graphs from any proper minor-closed graph class.
\end{abstract}

\footnotetext[2]{School of Electrical Engineering and Computer Science, University of Ottawa, Ottawa, Canada (\href{mailto:vida.dujmovic@uottawa.ca}{\nolinkurl{vida.dujmovic@uottawa.ca}}). Research supported by NSERC.}

\footnotetext[3]{D\'epartement d'Informatique, Universit\'e libre de Bruxelles, Brussels, Belgium (\href{mailto:robert.hickingbotham@ulb.be}{\nolinkurl{robert.hickingbotham@ulb.be}}).}

\footnotetext[4]{School of Computer Science, Carleton University, Ottawa, Canada (\href{mailto:morin@scs.carleton.ca}{\nolinkurl{morin@scs.carleton.ca}}). Research supported by NSERC.}

\footnotetext[5]{School of Mathematics, Monash University, Melbourne, Australia (\href{mailto:david.wood@monash.edu}{\nolinkurl{david.wood@monash.edu}}). Research supported by the Australian Research Council and by NSERC.}

\renewcommand{\thefootnote}{\arabic{footnote}}
\section{Introduction}

Graph product structure theory describes graphs in complicated graph classes as subgraphs of products of simpler graphs (typically with bounded treewidth)\footnote{Let $\mathdefn{\NN}:=\{0,1,2,\dots\}$ and $\mathdefn{\NN^+}:=\NN\setminus\{0\}$. We consider finite simple undirected graphs $G$ with vertex-set $V(G)$ and edge-set $E(G)$. For a vertex $v\in V(G)$, let \defn{$N_G(v)$} denote the set of neighbours of $v$ in $G$, and for $S\subseteq V(G)$, let $\mathdefn{N_G(S)}:=(\bigcup_{v\in S}N_G(v))\setminus S$ and $\mathdefn{N_G[S]}:=N_G(S)\cup S$. A graph $H$ is a \defn{minor} of a graph $G$ if a graph isomorphic to $H$ can be obtained from $G$ by a series of edge-deletions, vertex-deletions, and edge-contractions. A graph $G$ \defn{contains} a graph $H$ if some subgraph of $G$ is isomorphic to $H$. A graph $H$ is a \defn{topological minor} of a graph $G$ if $G$ contains a subdivision of $H$. A collection $\GG$ of graphs is a \defn{class} if it is closed under isomorphism. A class $\GG$ is \defn{minor-closed} if for every graph $G\in\GG$ every minor of $G$ is in $\GG$. The \defn{strong product} of two graphs $A$ and $B$, denoted by \defn{$A\boxtimes B$}, is the graph with vertex-set $V(A)\times V(B)$, where distinct vertices $(v,x)$ and $(w,y)$ are adjacent if $v=w$ and $xy\in E(B)$, or $vw\in E(A)$ and $x=y$, or $vw\in E(A)$ and $xy\in E(B)$.}. The classical result in this direction is the \defn{Planar Graph Product Structure Theorem} by \citet{DJMMUW20}:

\begin{thm}[{\protect Planar Graph Product Structure Theorem~\citep{DJMMUW20,UWY22}}]
\label{PGPST}
Every planar graph is contained in:
\begin{itemize}
    \item $H\boxtimes P$ for some planar graph $H$ with treewidth 6 and some path $P$,
    \item $H\boxtimes P\boxtimes K_2$ for some planar graph $H$ with treewidth 4 and some path $P$, 
    \item $H\boxtimes P\boxtimes K_3$ for some planar graph $H$ with treewidth 3 and some path $P$.
\end{itemize}
\end{thm}

This result has been the key tool in resolving several longstanding open problems~\cite{DJMMUW20,DEJWW20,DEGJMM21,BKW26,DEMW23,DGLTU22,DFMS21}. In most of these applications, the best bounds are obtained by applying the third version of the Planar Graph Product Structure Theorem, since the dependence on the treewidth of $H$ is often exponential.\footnote{Queue-number is a good example of this phenomenon. The best known upper bound on the queue-number of $H\boxtimes P \boxtimes K_c$ in terms of $\tw(H)$ is $3c\,2^{\tw(H)}+\floor{3c/2}$; see \citep{DJMMUW20}.} The treewidth 3 bound in this version is best possible since \citet{DJMMUW20} showed that for any integer $c\geq 1$ there is a planar graph $G$ such that for any graph $H$ and path $P$, if $G$ is contained in $H\boxtimes P\boxtimes K_c$, then $K_4\subseteq H$ implying $\tw(H)\geq 3$.

Graph product structure theorems analogous to \cref{PGPST} have been established for numerous graph classes of interest \cite{DJMMUW20,DHHW22,DMW23,HW24,DHSW24,ISW24,BDHK24,HJ24,HJ24a,LNW,HWY25,HKW,HW25a,DHKWY25}. Here we focus on minor-closed graph classes. \citet{DJMMUW20} proved the following result when an apex graph\footnote{A graph $X$ is \defn{apex} if $V(X)=\emptyset$ or $X-v$ is planar for some vertex $v\in V(X)$.} is excluded as a minor:

\begin{thm}[{\protect\citep[Theorem~4]{DJMMUW20}}]
\label{apex_minor_free}
    For every apex graph $X$ there exists $c\in\NN$ such that every $X$-minor-free graph $G$ is contained in $H\boxtimes P$ for some graph $H$ of treewidth at most $c$ and some path $P$.
\end{thm}

This result qualitatively generalises the Planar Graph Product Structure Theorem since $K_5$ (or $K_{3,3}$) is apex and non-planar.
In fact, \citet{DJMMUW20} showed that graphs in a fixed minor-closed graph class $\GG$ are contained in $H\boxtimes P\boxtimes K_c$ for some $H$ with bounded treewidth if and only if some apex graph is not in $\GG$ (see \citep{ISW24,DHHJLMMRW25} for extensions). 
\citet{DEMWW22} showed that the condition that $X$ is an apex graph in \cref{apex_minor_free} can be dropped if we only consider $X$-minor-free graphs of bounded maximum degree:

\begin{thm}[{\protect\citep[Lemma~15]{DEMWW22}}]
\label{old_main}
  For every graph $X$ and every $\Delta\in\NN$, there exists $c\in \NN$ such that every $X$-minor-free graph $G$ of maximum degree at most $\Delta$ is contained in $H\boxtimes P$ for some graph $H$ of treewidth at most $c$ and some path $P$.
\end{thm}

As discussed above, the primary quantitative goal in any $H\boxtimes P \boxtimes K_c$ product structure theorem is to minimise the treewidth of $H$. The main result in this paper is the following variant of \cref{old_main}, which makes the treewidth of $H$ independent of $X$ and $\Delta$:

\begin{thm}\label{main}
    For every graph $X$ and every $\Delta\in\NN$, there exists $c\in\NN^+$ such that every $X$-minor-free graph $G$ of maximum degree at most $\Delta$ is contained in $H\boxtimes P\boxtimes K_c$ for some graph $H$ of treewidth at most $3$ and some path $P$.    
\end{thm}

\textbf{Application to Weak Colouring Numbers}

Weak colouring numbers are an important family of graph parameters in structural and algorithmic graph theory; see the recent survey by \citet{Siebertz26}. Let $G$ be a graph and let $\sigma$ be an ordering of $V(G)$. 
For an integer $r\in \NN$ and two vertices $u$ and $v$ of $G$, we say that
$u$ is \defn{weakly $r$-reachable} from $v$ in $\sigma$, if there exists
a $uv$-path of length at most $r$ such that for every vertex $w$ on the path, $u\leq_{\sigma} w$.
Let \defn{$\WReach_r[G, \sigma, v]$} denote the set of vertices that are weakly $r$-reachable from a vertex $v$ in $\sigma$. The \defn{$r$-weak colouring number $\wcol_r(G)$} of $G$ is
\[
\wcol_r(G) = \min_{\sigma}\max_{v \in V(G)}\ |\WReach_r[G, \sigma, v]|,
\]
where $\sigma$ ranges over the set of all vertex orderings of $G$.
Several papers give bounds for weak colouring numbers of graphs in a given sparse graph class. For example, if $G$ is planar, then $\wcol_r(G)=O(r^3)$; 
if $G$ has no $K_t$-minor, then $\wcol_r(G)=O(r^{t-1})$ as proved by \citet{HOQRS17}; and if $\tw(G)\leq t$, then 
$\wcol_r(G)\leq \binom{r+t}{t}$ and this bound is exactly tight as proved by \citet{GKRSS18}. More recently, \citet{HLMR26b} determined, up to an $O(r)$ factor, the maximum $r$-weak colouring number of a graph within any minor-closed class.

The following lemma is obtained by combining Lemma~22 of Hickingbotham and Wood~\cite{HW24} with the above result of Grohe, Kreutzer, Rabinovich, Siebertz, and Stavropoulos~\cite{GKRSS18}.

\begin{lem}[\cite{HW24}]
     For each $t\in\NN$ and $c,r\in\NN^+$, if $H$ is a graph of treewidth at most $t$ and $P$ is a path, then $\wcol_r(H\boxtimes P\boxtimes K_c)=O_{t,c}(r^{t+1})$.
\end{lem}

Since $\wcol_r(G)$ is monotone, combining this lemma with \cref{main} gives the following application. In particular, the exponent in the polynomial is independent of the excluded minor. We are not aware of any polynomial bounds with this property for bounded-degree minor-free graphs.

\begin{thm}\label{WeakColouringApplcation}
    For each $\Delta\in\NN$ and $t,r\in\NN^+$, every $K_t$-minor-free graph $G$ of maximum degree at most $\Delta$ has $\wcol_r(G)=O_{t,\Delta}(r^{4})$.    
\end{thm}

\subsection{Proof Outline}

We now give a brief outline of the proof of \cref{main}.  
Say that a graph class $\calG$ has \defn{product structure} if there exists an integer $t_\calG$ such that, for each $G\in\calG$, there exists a graph $H$ of treewidth at most $t_\calG$ and a path $P$ such that $G$ is isomorphic to a subgraph of $H\boxtimes P$.   
We begin by explaining how \cref{old_main} has its roots in \cref{apex_minor_free}, which states that, for every apex graph $X$, the class of $X$-minor-free graphs has product structure.  The proof of \cref{apex_minor_free} by \citet{DJMMUW20} makes use of (a version of) the \emph{Graph Minor Structure Theorem}, which describes any graph excluding a fixed minor $X$ in terms of a tree-decomposition whose torsos are ``almost-embedded.''  Very roughly, each almost-embedded torso $Q$ consists of a graph $Q_0$ embedded on a surface of small Euler genus, a small number of non-embedded subgraphs $Q_1,\ldots,Q_p$ (called vortices), a small number of \defn{minor apex vertices} that are adjacent only to vertices of $Q_1,\ldots,Q_p$, and a small number of \defn{major apex vertices} that can be adjacent to any vertex of $Q$.  Here, ``small number'' means a quantity depending only on the excluded graph $X$.  The problematic vertices are the major apex vertices; any one of these vertices may be adjacent to an arbitrary subset of vertices in the torso.  As a consequence, the presence of even a single major apex vertex in a torso precludes the possibility of a torso having product structure.

Therefore, to prove \cref{apex_minor_free}, \citet{DJMMUW20} use a variant of the Graph Minor Structure Theorem which states that, if $X$ is an apex graph, then the torsos in the tree-decomposition of an $X$-minor-free graph need not have major apex vertices.  This allows them to show that each torso has product structure and allows them to combine torsos to show that the entire graph has product structure.

\Citet{DEMWW22} prove \cref{old_main} by using the fact that, in an $X$-minor-free graph of maximum degree at most $\Delta$, the tree-decomposition obtained from (an appropriate version of) the Graph Minor Structure Theorem can be modified so that all major apex vertices become minor apex vertices.  This increases the parameters of the embedding (the genus and/or number of vortices) but allows them to apply the results in \citet{DJMMUW20} for tree-decompositions whose torsos have no major apex vertices.  The end result is that, for any $\Delta$ and \emph{any} graph $X$, there exists an integer $t_{X,\Delta}$ such that each $X$-minor-free graph $G$ of maximum degree at most $\Delta$ is isomorphic to a subgraph of $H\boxtimes P$ for some graph $H$ of treewidth at most $t_{X,\Delta}$.

Our main result, \cref{main}, shows that the treewidth of $H$ can be reduced to $3$, independent of $X$ and $\Delta$, at the cost of a third factor $K_c$ in the product, where only $c$ depends on $X$ and $\Delta$.  The proof makes use of two main tools. The first is a result of \citet{UTW}:
\begin{lem}[{\protect\citep[Proposition~56~(arXiv~version)]{UTW}}]\label{CampbellProduct}
  For any $n,k\in\NN^+$, there exists $c\in O(n^3k^2)$ such that every $K_{3,n}$-topological-minor-free graph $G$ of treewidth at most $k$ is contained in $H\boxtimes K_c$ for some graph $H$ of treewidth at most $3$.
\end{lem}
Since any subdivision of $K_{3,n}$ contains vertices of degree $n$, a graph of maximum degree at most $\Delta$ cannot contain a subdivision of $K_{3,\Delta+1}$.  However, for $t\ge 5$, bounded-degree $K_t$-minor-free graphs have unbounded treewidth, so \cref{CampbellProduct} does not immediately apply to $G$. 

We can get closer to making use of \cref{CampbellProduct} by considering a formulation of product structure in terms of layered partitions:  A graph $G$ is isomorphic to a subgraph of $H\boxtimes P\boxtimes K_c$ if and only if there exists a vertex partition $\calP$ of $G$ and a layering $\calL$ of $G$ such that the quotient graph $G/\calP$ is isomorphic to a subgraph of $H$ and $|L\cap P|\le c$ for each $(L,P)\in\calL\times\calP$.   Thus, one might attempt to apply \cref{CampbellProduct} to the bounded-treewidth quotient graph $H=G/\calP$ implied by \cref{old_main}.  However, the quotient graph $H$ resulting from \cref{old_main} may not have bounded maximum degree even though $G$ has bounded maximum degree.\footnote{This is an unavoidable phenomenon: \citet{DJMMW24} showed that there is a class of planar graphs of maximum degree $5$ for which any product structure result in which the treewidth of $H$ is bounded requires that $H$ have unbounded maximum degree.} As a result, it is not immediately clear that $H$ excludes $K_{3,n}$-subdivisions for any bounded $n$.  Indeed, a quotient graph $H$ produced by an application of \cref{old_main} could easily contain a subdivision of $K_{3,n}$ for arbitrarily large values of $n$.  

To overcome this obstacle, we adapt the ``tree-of-curtains'' decomposition introduced by \citet{DEMW23} in their proof of the Clustered Hadwiger Conjecture.  We partition the tree-decomposition obtained from the Graph Minor Structure Theorem into subtrees (each defining a subgraph of $G$ called a \defn{curtain} in \cite{DEMW23}).  For each such subtree, we show that the corresponding curtain is isomorphic to a subgraph of $J\boxtimes P$, where $P$ is a path and $J$ is a graph of bounded treewidth that contains no subdivision of $K_{3,n}$ for some bounded $n$.  This allows us to apply \cref{CampbellProduct} to the quotient graph $J$ to conclude that each curtain is isomorphic to a subgraph of $H\boxtimes P\boxtimes K_c$, where $H$ has treewidth at most $3$ and $c$ is bounded.  

By carefully modifying the almost-embeddings of the torsos in the Graph Minor Structure Theorem before defining the tree-of-curtains, we can ensure that the adhesion between any pair of adjacent curtains consists of at most three vertices of $G$.  This allows us to set up an inductive proof in which a root curtain $G_0$ is adjacent to a collection of subgraphs $G_1,\ldots,G_s$. A sufficiently strong inductive hypothesis allows us to show that each subgraph $G_i$ is isomorphic to a subgraph of $H_i\boxtimes P_i\boxtimes K_c$, where $P_i$ is a path and $H_i$ is a graph of treewidth at most $3$.  Furthermore, the at most three vertices in the adhesion $S_i:=V(G_0)\cap V(G_i)$ form a single part in the partition that defines the quotient graph $H_i$.  We add edges to $G_0$ to obtain a graph $G^+_0$ in which the vertices in $S_i$ form a clique for each $i\in\{1,\ldots,s\}$.  Since $G^+_0$ is a curtain, it is isomorphic to a subgraph of $H_0\boxtimes P_0\boxtimes K_c$, where $H_0$ has treewidth at most $3$.  Furthermore, for each $i\in\{1,\ldots,s\}$, the vertices in $S_i$ correspond to at most $3$ vertices of $H_0$ that form a clique in $H_0$.  This allows us to merge $H_0,H_1,\ldots,H_s$ via clique-sums of size at most $4$ so that the treewidth of the resulting graph $H$ is at most $\max\{3,\max_{0\le i\le s} \tw(H_i)\} \le 3$.    

The remainder of this paper is organized as follows.
\Cref{tools} develops the tools we use, the most notable of which is \cref{K3nBound}, which bounds the maximum value of $n$ for which a curtain can contain a subdivision of $K_{3,n}$. \Cref{proof} then uses these tools to prove \cref{main}.  

\section{Tools}
\label{tools}

\subsection{Surface Lemmas}

A \defn{surface} is a compact connected 2-dimensional manifold without boundary. For any terms related to graphs on surfaces that are not defined here, the reader is referred to the monograph by \citet{MoharThom}.

Let $G$ be a graph embedded in a surface $\Sigma$. For distinct vertices $a,b\in V(G)$, two internally disjoint $ab$-paths $P_1$ and $P_2$ in $G$ are \defn{homotopic} if $P_1 \cup P_2$ bounds a disc in $\Sigma$. (If $\Sigma$ has Euler genus $0$, that is, $\Sigma$ is a sphere, then any two internally disjoint $ab$-paths are homotopic.) The following lemma shows that any sufficiently large set of pairwise internally-disjoint $ab$-paths contains a large set of pairwise homotopic $ab$-paths.

\begin{lem}[{\protect\citep[Lemma~17]{HW26}}]
\label{HomotopicTuple}
  For any integer $c\ge 2$, any graph $G$ embedded in a surface of Euler genus $g\ge 0$, and any distinct vertices $a,b\in V(G)$, if $\mathcal{P}$ is a set of pairwise internally disjoint $ab$-paths in $G$ and $|\mathcal{P}| \ge \max\{c, (c-1)(2g+1)\}$, then $\mathcal{P}$ contains a set of $c$ pairwise homotopic paths.
\end{lem}

Let $G$ be a graph embedded on a surface $\Sigma$ and let $\calP'$ be a set of pairwise homotopic $ab$-paths in $G$.  Let $D$ be a closed disc in $\Sigma$ whose boundary is the union of two paths in $\calP'$ and that maximizes the number of paths in $\calP'$ contained in $D$.  Then $D$ must contain all of the paths in $\calP'$ since any path in $\calP'$ whose interior is disjoint from $D$ could be used to construct a disc $D'\supseteq D$ that contradicts the maximality of $D$.  Applying this argument to the set $\calP'$ guaranteed by \cref{HomotopicTuple} leads to the following:

\begin{cor}
\label{homotopic_disc}
  For any integer $c\ge 2$, any graph $G$ embedded in a surface $\Sigma$ of Euler genus $g\ge 0$, and any distinct vertices $a,b\in V(G)$, if $\mathcal{P}$ is a set of pairwise internally disjoint $ab$-paths in $G$ and $|\mathcal{P}| \ge \max\{c, (c-1)(2g+1)\}$, then $\mathcal{P}$ contains a subset $\mathcal{P}'$ of $c$ paths that are all contained in a disc $D\subset\Sigma$ whose boundary is the union of two paths in $\calP'$.
\end{cor}

\subsection{Graph Minor Structure Theorem}

The Graph Minor Structure Theorem of \citet{RS-XVI} describes the structure of graphs excluding a fixed minor using four ingredients: tree-decompositions, graphs on surfaces, vortices, and apex vertices. To describe this formally, we need the following definitions.

A \defn{tree-decomposition} of a graph $G$ is a collection $(B_x :x\in V(T))$ of subsets of $V(G)$ (called \defn{bags}) indexed by the vertices of a tree $T$, such that (a) for every edge $uv\in E(G)$, some bag $B_x$ contains both $u$ and $v$, and (b) for every vertex $v\in V(G)$, the set $\{x\in V(T):v\in B_x\}$ induces a non-empty (connected) subtree of $T$. For each edge $xy\in E(T)$, the set $B_x\cap B_y$ is called an \defn{adhesion}. 
The \defn{width} of $(B_x:x\in V(T))$ is $\max\{|B_x| \colon x\in V(T)\}-1$. A \defn{path-decomposition} is a tree-decomposition in which the underlying tree is a path, simply denoted by the corresponding sequence of bags $(B_1,\dots,B_n)$. The \defn{treewidth} of a graph $G$, denoted by \defn{$\tw(G)$}, is the minimum width of a tree-decomposition of $G$.  

If $(B_x:x\in V(T))$ is a tree-decomposition of a graph $G$, then the \defn{torso} \defn{$\torso{G}{B_x}$} of a bag $B_x$ is the graph obtained from the induced subgraph $G[B_x]$ by adding edges so that $B_x\cap B_y$ is a clique for each edge $xy\in E(T)$.

A \defn{rooted tree} consists of a tree $T$ and a distinguished vertex of $T$ called the \defn{root}.  For a rooted tree $T$ and a node $x$ of $T$, let \defn{$T_x$} denote the subtree of $T$ induced by $x$ and all the descendants of $x$. The \defn{depth} of a vertex $v$ in a tree $T$ rooted at $r\in V(T)$ is $\dist_T(r,v)$. A tree-decomposition $(B_x:x\in V(T))$ of a graph $G$ is \defn{rooted} if $T$ is rooted. If $T$ is rooted at $r\in V(T)$, then $\torso{G}{B_r}$ is called the \defn{root torso}.  

Let $G_0$ be a graph embedded in a surface $\Sigma$. A closed topological disc $D$ in $\Sigma$ is \defn{$G_0$-clean} if the interior of $D$ is disjoint from $G_0$, and the boundary of $D$ only intersects $G_0$ in vertices of $V(G_0)$. Let $x_1,\dots,x_n$ be the vertices of $G_0$ on the boundary of $D$ in order around $D$. A \defn{$D$-vortex} consists of a graph $H$ and a path-decomposition $(B_1,\dots,B_n)$ of $H$ such that $x_i\in B_i$ for each $i\in\{1,\dots,n\}$, and $V(G_0\cap H)=\{x_1,\dots,x_n\}$. The vertices $x_1,\dots,x_n$ are called the \defn{boundary vertices} of the $D$-vortex. (When $n=0$, a $D$-vortex is the $0$-vertex graph $H$ with the empty path-decomposition $()$, which by convention has width $0$.)

For integers $g,p,k,a\geq 0$, a \defn{$(g,p,k,a)$-almost-embedding} of a graph $G$ is a tuple $\Phi:=(G_0,G_1,\dots,G_p,A)$ such that:
\begin{itemize}
  \item  $A\subseteq V(G)$ with $|A|\le a$;
  \item $G_0,G_1,\dots,G_p$ are subgraphs of $G$ such that $G-A = G_{0} \cup G_{1} \cup \cdots \cup G_p$;
  \item $G_{1}, \dots, G_p$ are pairwise vertex-disjoint;
  \item $G_{0}$ is embedded in a surface $\Sigma$ of Euler genus at most $g$;
  \item there are pairwise disjoint $G_0$-clean topological discs $D_1,\dots,D_p$ in $\Sigma$; and
  \item $G_i$ is a $D_i$-vortex with a path-decomposition $(B_1,\dots,B_{n_i})$ of width at most $k$, for each $i\in\{1,2,\dots,p\}$.
\end{itemize}

The vertices in $A$ are called \defn{apex} vertices, and \defn{$\apex(\Phi)$} denotes the set of apex vertices of an almost-embedding $\Phi$. The \defn{top} of $\Phi$ is $\mathdefn{\topp(\Phi)}:=A\cup V(G_1\cup \cdots\cup G_p)$.  
For a collection of graphs, each equipped with a $(g,p',k,a)$-almost-embedding with $p'\le p$, we use the convention that $p'=p$ for each of the almost-embeddings, by introducing $0$-vertex vortices.  (We do this only to avoid additional notation.)

An \defn{almost-embedded tree-decomposition} of a graph $G$ is a pair $(\calT,\calE)$ where $\calT:=(B_x:x\in V(T))$ is a tree-decomposition of $G$, and $\calE:=(\Phi_x:x\in V(T))$ where $\Phi_x$ is an almost-embedding of the torso $\torso{G}{B_x}$ for each $x \in V(T)$. If each $\Phi_x$ is a $(g,p,k,a)$-almost-embedding, then $(\calT,\calE)$ is a \defn{$(g,p,k,a)$-almost-embedded tree-decomposition}. We make use of the following variant of the Graph Minor Structure Theorem, due to \citet{DKMW12}, which focuses on the structure of the adhesions.  

\begin{thm}[{\protect\citep[Theorem~4]{DKMW12}}]
\label{GMST-DKMW}
For every integer $t\ge 1$, there exist integers $g,p,k,a\geq 0$ such that every $K_t$-minor-free graph $G$ has a rooted $(g,p,k,a)$-almost-embedded tree-decomposition $(\calT,\calE):=((B_x:x\in V(T)), (\Phi_x:x\in V(T)))$ such that for each node $y$ of $T$ with parent $x$, $B_x\cap B_y\subseteq\apex(\Phi_y)$ and $|(B_x\cap B_y)\setminus\topp(\Phi_x)|\le 3$.
\end{thm}

Indeed, in \citet[Theorem~4]{DKMW12}, for each node $y$ with parent $x$, $B_x\cap B_y\subseteq\apex(\Phi_y)$, and either $(B_x\cap B_y)\setminus\apex(\Phi_x)$
is contained in a single bag of the path-decomposition of some vortex of $\Phi_x$ (which implies $B_x\cap B_y\subseteq\topp(\Phi_x)$, so $|(B_x\cap B_y)\setminus\topp(\Phi_x)|=0$), or $(B_x\cap B_y)\setminus\apex(\Phi_x)$ consists of at most $3$ vertices lying on a single face of the embedded subgraph of $\Phi_x$ (which implies $|(B_x\cap B_y)\setminus\topp(\Phi_x)|\le |(B_x\cap B_y)\setminus\apex(\Phi_x)|\le 3$).
 
On the way to proving \cref{main}, the almost-embedded tree-decomposition given by \cref{GMST-DKMW} will go through several stages of cleaning up.  The rest of this subsection shows how to implement these.

The following shows that each adhesion $B_x\cap B_y$ in a $(g,p,k,a)$-almost-embedded tree-decomposition has size at most $\mks$.\footnote{The bound in \citep[Lemma~21]{DMW17} is actually $a + 2k + \floor{ \frac12 (7 + \sqrt{1+24g})}$ which is less than $8\max\{2,g,k,a\}$ for all $g,k,a$.} 

\begin{lem}[{\protect\citep[Lemma~21]{DMW17}}]
\label{mks}
  If a graph $G$ has a $(g,p,k,a)$-almost-embedding, then every clique in $G$ has at most $\mks$ vertices.
\end{lem}

A rooted tree-decomposition $(B_x:x\in V(T))$ is \defn{crisp} if, for each node $x$ of $T$, each child $y$ of $x$, and each vertex $v$ in the adhesion $B_x\cap B_y$, there exists an edge $vw$ of $G$ with $w\in (\bigcup_{z\in V(T_y)} B_z)\setminus B_x$.  Intuitively, there is no reason for a tree-decomposition to be non-crisp, since there is no reason for the vertex $v$ to be present in $B_z$ for any $z\in V(T_y)$.

\begin{lem}\label{crisp_degree}
  Let $G$ be a graph, and let $\calT=(B_x:x\in V(T))$ be a crisp rooted tree-decomposition of $G$.  Then, for each $v\in V(G)$ and each $x\in V(T)$, the number of children $y$ of $x$ such that $v\in B_x\cap B_y$ is at most $\deg_G(v)$.  
\end{lem}

\begin{proof}
  Let $x\in V(T)$ and let $v\in V(G)$.  Let $Y$ be the set of children $y$ of $x$ such that $v\in B_x\cap B_y$.  For each $y\in Y$, $G$ contains an edge $vw_y$ with $w_y\in (\bigcup_{z\in V(T_{y})} B_z)\setminus B_x$.  Since the sets in $\{(\bigcup_{z\in V(T_{y})} B_z)\setminus B_x : y\in Y\}$ are pairwise disjoint, $\deg_G(v)\ge |\{w_y:y\in Y\}|= |Y|$, as required.
\end{proof}

The following lemma shows that a rooted almost-embedded tree-decomposition whose parent adhesions are contained in apex sets can be made crisp.

\begin{lem}\label{chordal_closure}
  Let $g,p,k,a\in\NN$, let $G$ be a graph, and let $(\calT_0,\calE_0):=((B^0_x:x\in V(T)),(\Phi^0_x:x\in V(T)))$ be a rooted $(g,p,k,a)$-almost-embedded tree-decomposition of $G$ such that for each node $y$ of $T$ with parent $x$, $B^0_x\cap B^0_y\subseteq\apex(\Phi^0_y)$.  Then $G$ has a crisp rooted $(g,p,k,a)$-almost-embedded tree-decomposition $(\calT,\calE):=((B_x:x\in V(T)),(\Phi_x:x\in V(T)))$ such that for each node $y\in V(T)$,
  \begin{enumerate}[label=\rm(\roman*)]
    \item\label{closure_bag_subset} $B_y \subseteq B^0_y$;
    \item\label{closure_apex_subset} $\apex(\Phi_y) = \apex(\Phi^0_y)\cap B_y$ and $\topp(\Phi_y) = \topp(\Phi^0_y)\cap B_y$; and
    \item\label{closure_adhesion_apex} if $y$ has parent $x$, then $B_x\cap B_y\subseteq\apex(\Phi_y)$ and $(B_x\cap B_y)\setminus\topp(\Phi_x)\subseteq (B^0_x\cap B^0_y)\setminus\topp(\Phi^0_x)$.
  \end{enumerate}
\end{lem}

\begin{proof} 
  We first define $(B_y:y\in V(T))$ by processing the nodes of $T$ in a bottom-up manner so that when we process a node $y$, each proper descendant $z\in V(T_y-y)$ has already been processed and satisfies $B_z\subseteq B^0_z$ as well as the following \defn{bottom-up invariant}: for each $u\in B^0_{p(z)}\cap B_z$ (where $p(z)$ is the parent of $z$), there exists $uw\in E(G)$ with $w\in (\bigcup_{q\in V(T_z)} B_q)\setminus B^0_{p(z)}$.

  If $y$ is the root of $T$, we set $Z_y:=\emptyset$ and $B_y:=B^0_y$. Otherwise, let $x$ be the parent of $y$, and let
  \[
    Z_y := \{v\in B^0_x\cap B^0_y:(N_G(v)\setminus B^0_x)\cap (B^0_y \cup \textstyle\bigcup_{z\in V(T_y-y)} B_z)=\emptyset\} \enspace ,
  \]
  and set $B_y:=B^0_y\setminus Z_y$.  Note that $Z_y\subseteq\apex(\Phi^0_y)$ for each $y\in V(T)$.  We claim that $Z_y\cap B_z=\emptyset$ for every $z\in V(T_y-y)$.  Suppose not, and choose $z\in V(T_y-y)$ of minimum depth such that $Z_y\cap B_z$ is non-empty. Let $v\in Z_y\cap B_z$ and let $p$ be the parent of $z$. Since $v\in Z_y\subseteq B^0_y$ and $v\in B_z\subseteq B^0_z$, we have  $v\in B^0_p\cap B_z$. Hence there exists $vw\in E(G)$ with $w\in (\bigcup_{q\in V(T_z)} B_q)\setminus B^0_p$. Since $w\in \bigcup_{q\in V(T_z)} B^0_q$ and $w\notin B^0_p$, we have $w\notin B^0_x$, contradicting $v\in Z_y$. Thus $Z_y\cap B_z=\emptyset$ for all $z\in V(T_y-y)$. Moreover, every edge $uv\in E(G)$ with $v\in Z_y$ and $u\in B^0_y$ satisfies $u\in B^0_x$ by the definition of $Z_y$, so $\{u,v\}\subseteq B^0_x$. 
  Furthermore, for each $v\in B^0_x\cap B_y = (B^0_x\cap B^0_y)\setminus Z_y$, there exists $vw\in E(G)$ with $w\in (\bigcup_{z\in V(T_y)} B_z)\setminus B^0_x$, maintaining the bottom-up invariant.

  After processing all nodes of $T$, $\calT:=(B_y:y\in V(T))$ is a tree-decomposition of $G$ satisfying \ref{closure_bag_subset}. For each $y\in V(T)$, let $\Phi^0_y:=(G^y_0,G^y_1,\ldots,G^y_p,A^y)$. Since $Z_y\subseteq\apex(\Phi^0_y)=A^y$ (including when $y$ is the root, where $Z_y=\emptyset$), $Z_y$ is disjoint from $V(G^y_0\cup G^y_1\cup\cdots\cup G^y_p)=B^0_y\setminus A^y$. Since $B_z\subseteq B^0_z$ for all $z\in V(T)$, the torso $\torso{G}{B_y}$ is a spanning subgraph of $\torso{G}{B^0_y}-Z_y$. Setting $G'_i:=G^y_i\cap\torso{G}{B_y}$ for each $i\in\{0,1,\ldots,p\}$ preserves the vertex set $V(G'_i)=V(G^y_i)$, the surface embedding of $G'_0\subseteq G^y_0$, the $G'_0$-clean vortex discs, and the vortex path-decompositions of width at most $k$. Thus $\Phi_y:=(G'_0,G'_1,\ldots,G'_p,A^y\setminus Z_y)$ is a $(g,p,k,a)$-almost-embedding of $\torso{G}{B_y}$ with $\apex(\Phi_y)=A^y\setminus Z_y = \apex(\Phi^0_y)\cap B_y$ and $\topp(\Phi_y)=\topp(\Phi^0_y)\cap B_y$, satisfying \ref{closure_apex_subset}.

  For each non-root node $y$ with parent $x$, since $B_x\subseteq B^0_x$, we have $B_x\cap B_y\subseteq B^0_x\cap B_y = (B^0_x\cap B^0_y)\setminus Z_y\subseteq \apex(\Phi^0_y)\setminus Z_y = \apex(\Phi_y)$.  For each $v\in B_x\cap B_y$, the vertex $w\in (\bigcup_{z\in V(T_y)} B_z)\setminus B^0_x$ from the construction satisfies $w\in (\bigcup_{z\in V(T_y)} B_z)\setminus B_x$ because $B_x\subseteq B^0_x$; hence $\calT$ is crisp.  Finally, if $v\in B_x\cap B_y$ and $v\in\topp(\Phi^0_x)$, then since $v\in B_x$ we have $v\in\topp(\Phi^0_x)\cap B_x = \topp(\Phi_x)$. Therefore, $(B_x\cap B_y)\setminus\topp(\Phi_x)\subseteq (B^0_x\cap B^0_y)\setminus\topp(\Phi^0_x)$, satisfying \ref{closure_adhesion_apex}.
\end{proof}

Applying \cref{chordal_closure} to the tree-decomposition guaranteed by \cref{GMST-DKMW}, we obtain the following corollary:

\begin{cor}\label{dkmw_crisp}
  For every integer $t\ge 1$, there exist integers $g,p,k,a\ge 0$ such that each $K_t$-minor-free graph $G$ has a crisp rooted $(g,p,k,a)$-almost-embedded tree-decomposition $(\calT,\calE):=((B_x:x\in V(T)), (\Phi_x:x\in V(T)))$ such that for each node $y$ of $T$ with parent $x$, $B_x\cap B_y\subseteq\apex(\Phi_y)$ and $|(B_x\cap B_y)\setminus\topp(\Phi_x)|\le 3$.
\end{cor}

One final tool for manipulating almost-embedded tree-decompositions is the following, which allows us to avoid adhesions that mix apex and non-top vertices.

\begin{lem}\label{promote_apex_neighbours}
  Let $g,p,k,a,d\ge 0$ be integers, let $G$ be a graph, and let $(\calT,\calE):=((B_x:x\in V(T)),(\Phi_x:x\in V(T)))$ be a $(g,p,k,a)$-almost-embedded tree-decomposition of $G$ such that for each $x\in V(T)$ and each $u\in\apex(\Phi_x)$,
  \[
    |N_{\torso{G}{B_x}}(u)\setminus\topp(\Phi_x)|\le d \enspace .
  \]
  Then there exists a $(g+2ad,p+1,k,a)$-almost-embedded tree-decomposition $(\calT,\calF:=(\Psi_x:x\in V(T)))$ such that
  \begin{enumerate}[label=\rm(\roman*)] 
    \item\label{bigger_tops} $\topp(\Phi_x)\subseteq\topp(\Psi_x)$ for each $x\in V(T)$;
    \item\label{same_apices} $\apex(\Phi_x)=\apex(\Psi_x)$ for each $x\in V(T)$;
    \item\label{apex_neighbourhood} $N_{\torso{G}{B_x}}(\apex(\Psi_x))\subseteq\topp(\Psi_x)$ for each $x\in V(T)$.
  \end{enumerate}
\end{lem}

\begin{proof}
  For each $x\in V(T)$, let $(G^x_0,G^x_1,\ldots,G^x_p,A^x):=\Phi_x$ and let $\Sigma^x$ be the surface of Euler genus at most $g$ on which $G^x_0$ is embedded. For each $i\in\{1,\ldots,p\}$, let $D^x_i$ be the $G^x_0$-clean disc whose boundary contains the vertices in $V(G^x_i)\cap V(G^x_0)$. Let $N^x_0:=N_{\torso{G}{B_x}}(A^x)\setminus\topp(\Phi_x)$.  Observe that $|N^x_0|\le \sum_{u\in A^x}|N_{\torso{G}{B_x}}(u)\setminus\topp(\Phi_x)|\le ad$.
  
  If $N^x_0=\emptyset$, let $\Pi^x:=\Sigma^x$, let $D^x_{p+1}$ be a closed disc in $\Sigma^x\setminus (G^x_0\cup D^x_1\cup\cdots\cup D^x_p)$ (which is a non-empty open subset of $\Sigma^x$), and let $G^x_{p+1}$ be the $0$-vertex graph with the empty path-decomposition $()$ of width $0\le k$. Then $\Psi_x:=(G^x_0,G^x_1,\ldots,G^x_{p+1},A^x)$ is a $(g+2ad,p+1,k,a)$-almost-embedding of $\torso{G}{B_x}$.

  Now assume that $N^x_0\neq\emptyset$. Since each $v\in N^x_0$ belongs to $V(G^x_0)\setminus\topp(\Phi_x)$, $v$ does not participate in any vortex of $\Phi_x$ and is therefore disjoint from $D^x_1,\ldots,D^x_p$. Thus we can choose pairwise disjoint closed discs $(D_v:v\in N^x_0)$ in $\Sigma^x$, all disjoint from $D^x_1,\ldots,D^x_p$, such that $v\in\partial D_v$ and $D_v\cap G^x_0=\{v\}$ for each $v\in N^x_0$.
  Choose an arbitrary $v\in N^x_0$. Choose pairwise disjoint closed discs $(O_w:w\in N^x_0\setminus\{v\})$ in $\operatorname{int}(D_v)$, and for each $w\in N^x_0\setminus\{v\}$ choose a closed disc $O'_w$ in $\operatorname{int}(D_w)$ and attach a handle (cylinder) $C_w$ to $\Sigma^x$ by removing $\operatorname{int}(O_w)\cup\operatorname{int}(O'_w)$ and identifying the boundary circles of $C_w$ with $\partial O_w$ and $\partial O'_w$. The resulting surface $\Pi^x$ has Euler genus at most $g+2(|N^x_0|-1)\le g+2ad$, and $G^x_0$ remains embedded in $\Pi^x$. Let
  \[
    \Gamma^x := \Bigl(D_v\setminus\bigcup_{w\in N^x_0\setminus\{v\}}\operatorname{int}(O_w)\Bigr) \cup \bigcup_{w\in N^x_0\setminus\{v\}} \bigl(C_w\cup (D_w\setminus\operatorname{int}(O'_w))\bigr) \enspace .
  \]
  Then $\Gamma^x$ is connected, disjoint from $D^x_1,\ldots,D^x_p$, and satisfies $\Gamma^x\cap G^x_0=N^x_0\subseteq\partial\Gamma^x$. Taking a regular neighbourhood in $\Gamma^x$ of a star of $|N^x_0|-1$ pairwise internally disjoint arcs in $\operatorname{int}(\Gamma^x)\cup N^x_0$ joining $v$ to each $w\in N^x_0\setminus\{v\}$ yields a $G^x_0$-clean closed disc $D^x_{p+1}\subseteq\Gamma^x$ with $\partial D^x_{p+1}\cap G^x_0=N^x_0$. Let $x_1,\ldots,x_{|N^x_0|}$ be the vertices of $N^x_0$ in order around $\partial D^x_{p+1}$, and define $G^x_{p+1}$ to be the edgeless graph $(N^x_0,\emptyset)$, which is a $D^x_{p+1}$-vortex with path-decomposition $(\{x_1\},\ldots,\{x_{|N^x_0|}\})$ of width $0\le k$.  Then $\Psi_x:=(G^x_0,G^x_1,\ldots,G^x_{p+1},A^x)$ is a $(g+2ad,p+1,k,a)$-almost-embedding of $\torso{G}{B_x}$.  
  
  The pair $(\calT,\calF:=(\Psi_x:x\in V(T)))$ satisfies \ref{bigger_tops} because $\topp(\Psi_x)=\topp(\Phi_x)\cup N^x_0$, satisfies \ref{same_apices} since $\apex(\Psi_x)=A^x=\apex(\Phi_x)$ for each $x\in V(T)$, and satisfies \ref{apex_neighbourhood} since $N_{\torso{G}{B_x}}(\apex(\Psi_x))=N_{\torso{G}{B_x}}(A^x)\subseteq\topp(\Phi_x)\cup N^x_0 = \topp(\Psi_x)$ for each $x\in V(T)$.
\end{proof}

\subsection{Layered Partitions}

A \defn{layering} $\calL:=(L_1,L_2,\ldots)$ of a graph $G$ is a sequence of pairwise disjoint, possibly empty subsets of $V(G)$ whose union is $V(G)$ and such that, for each edge $vw\in E(G)$, $v\in L_i$ and $w\in L_j$ implies that $|i-j|\le 1$.  For a graph $H$ and a graph $G$, an \defn{$H$-partition} $\calP:=(C_x:x\in V(H))$ of $G$ is a collection of pairwise disjoint, possibly empty subsets of $V(G)$ whose union is $V(G)$ such that, for each edge $vw\in E(G)$, $v\in C_x$ and $w\in C_y$ implies that $x=y$ or $xy\in E(H)$.  The \defn{width} of $\calP$ is $\max(\{0\}\cup\{|C_x|:x\in V(H)\})$.  (Note that a layering of $G$ is an $H$-partition of $G$ where $H:=(1,2,\ldots)$ is a path graph.) 

An $H$-partition $(C_x:x\in V(H))$ is \defn{connected} if $C_y$ is non-empty and $G[C_y]$ is connected for each $y\in V(H)$.  An $H$-partition $(C_x:x\in V(H))$ is \defn{minimal} if it has no empty parts and, for each edge $xy\in E(H)$, there exists $vw\in E(G)$ with $v\in C_x$ and $w\in C_y$.\footnote{When an $H$-partition $\calP:=(C_x:x\in V(H))$ of a graph $G$ is minimal, $H$ is isomorphic to the quotient graph $G/\calP$.} For convenience, we do not distinguish between the vertex $v$ in a singleton part $\{v\}\in \calP$ and the node $y\in V(H)$ with $C_y=\{v\}$ that indexes $\{v\}$.  In other words, $y=v$ whenever $C_y=\{v\}$.  We make use of the following easy observation:

\begin{obs}\label{singleton_degree}
  Let $\calP:=(C_x:x\in V(H))$ be a minimal connected $H$-partition of a graph $G$. Then, for each singleton part $C_v=\{v\}\in\calP$, $\deg_H(v)\le\deg_G(v)$.  
\end{obs}

A \defn{layered $H$-partition} of a graph $G$ is a pair $(\calL, \calP)$ where $\calL:=(L_1,L_2,\ldots)$ is a layering of $G$, and $\calP:=(C_y:y\in V(H))$ is an $H$-partition of $G$.  The \defn{layered width} of $(\calL,\calP)$ is $\max(\{0\}\cup\{|L_i\cap C_y|:i\in\NN^+,\, y\in V(H)\})$.  A layered $H$-partition $(\calL,\calP)$ of a graph $G$ is \defn{connected} (respectively, \defn{minimal}) if $\calP$ is connected (respectively, minimal).  Layered partitions are closely related to \cref{main}: A graph $G$ has a layered $H$-partition of layered width at most $c$ if and only if $G$ is isomorphic to a subgraph of $H\boxtimes P\boxtimes K_{c}$, where $P$ is a path.

We make use of the following result, a version of which appears in \cite{DEMW23} and follows quickly from \cite[Theorem~38]{DJMMUW20}:

\begin{lem}[\citep{DEMW23,DJMMUW20}]
\label{almost_embedded_layered_partition}
  There exists a function $\ell:\NN^4\to\NN^+$ such that the following holds, for each $g,p,k,a\in\NN$.  Let $\Phi$ be a $(g,p,k,a)$-almost-embedding of a graph $G$ with the property that $N_G(\apex(\Phi))\subseteq\topp(\Phi)$.  Then $G$ has a minimal connected layered $H$-partition $(\calL,\calP):=((L_1,L_2,\ldots), (C_x:x\in V(H)))$ of layered width at most $\ell(g,p,k,a)$, where the treewidth of $H$ is at most $\ell(g,p,k,a)$, and such that
  \begin{enumerate}[label=\rm(\roman*),nosep,nolistsep]
    \item $\topp(\Phi)\subseteq L_1$;
    \item for each $v\in L_1\cup L_2$, the part $C_{v}=\{v\}$ of $\calP$ that contains $v$ is a singleton.
  \end{enumerate}
\end{lem}

It will be helpful to work with a weaker notion of layering and a corresponding weaker notion of layered partition.  For $c\in\NN^+$, a \defn{$c$-weak layering} $(L_1,L_2,\ldots)$ of a graph $G$ is a sequence of pairwise disjoint subsets of $V(G)$ that cover $V(G)$ such that for each $vw\in E(G)$, $v\in L_i$ and $w\in L_j$ implies $|i-j|\le c$.  A \defn{$c$-weak layered $H$-partition} of a graph $G$ is a pair $(\calL, \calP)$ where $\calL:=(L_1,L_2,\ldots)$ is a $c$-weak layering of $G$, and $\calP:=(C_y:y\in V(H))$ is an $H$-partition of $G$.  The \defn{layered width} of $(\calL,\calP)$ is $\max(\{0\}\cup\{|L_i\cap C_y|:i\in\NN^+,\, y\in V(H)\})$.  The following lemma shows that $c$-weak layered partitions can be converted to layered partitions with a factor $c$ increase in the layered width.

\begin{lem}\label{weak_to_strong}
  If a graph $G$ has a $c$-weak layered $H$-partition of layered width at most $w$, then $G$ has a layered $H$-partition of layered width at most $cw$.
\end{lem}

\begin{proof}
  Let $(\calL:=(L_1,L_2,\ldots),\calP)$ be a $c$-weak layered $H$-partition of layered width at most $w$. For each $i\in\NN^+$, let $L'_i:=\bigcup_{j=(i-1)c+1}^{ic} L_j$. Then $\calL':=(L'_1,L'_2,\ldots)$ is a sequence of pairwise disjoint sets covering $V(G)$.  Furthermore, for each edge $vw$ of $G$, $v\in L_x$ and $w\in L_y$ implies that $|x-y|\le c$. If $v\in L'_i$ and $w\in L'_j$, then $x \in \{(i-1)c+1, \ldots, ic\}$ and $y \in \{(j-1)c+1, \ldots, jc\}$. Since $|x-y| \le c$, it follows that $|i-j| \le 1$, meaning $\calL'$ is a proper layering of $G$. For each $i\in\NN^+$ and $y\in V(H)$, we have $|L'_i\cap C_y|\le \sum_{j=(i-1)c+1}^{ic} |L_j\cap C_y|\le cw$.  Therefore $(\calL',\calP)$ is a layered $H$-partition of layered width at most $cw$.
\end{proof}

\subsection{Graphs without \texorpdfstring{$K_{3,n}$}{K\_\{3,n\}} Subdivisions}

In this subsection, some of the heavy lifting begins. We will establish a lemma that bounds the maximum size of a $K_{3,n}$ topological minor in a graph with appropriately constrained almost-embedded tree-decompositions.  Our proof of \cref{main} will use this lemma so that we can apply \cref{CampbellProduct} (the result of \citet{UTW} giving $H$-partitions of bounded-treewidth graphs having no $K_{3,n}$-subdivision) to specific subgraphs of the $K_t$-minor-free graphs considered by \cref{main}.  Each such subgraph will be defined by a subtree of an almost-embedded tree-decomposition $(\calT,\calE)$ obtained by applying \cref{promote_apex_neighbours} to an almost-embedded tree-decomposition of a $K_t$-minor-free graph obtained from \cref{dkmw_crisp}.

\begin{lem}
\label{MonoK2n}
  For every integer $n\ge 1$, every edge-2-colouring of $K_{3,n}$ has a monochromatic $K_{2,\lceil n/8\rceil}$.
\end{lem}

\begin{proof}
  Let $(\{v_1,v_2,v_3\},\{w_1,\dots,w_{n}\})$ be the bipartition of $K_{3,n}$. For each $i\in\{1,\ldots,n\}$, label $w_i$ by 
  \[
    ( \col(v_1w_i), \col(v_2w_i), \col(v_3w_i) ) \enspace .
  \]
  There are eight possible labels. By the Pigeonhole Principle, some label appears at least $\lceil n/8\rceil$ times. Let $W$ be a set of $\lceil n/8\rceil$ vertices in $\{w_1,\dots,w_{n}\}$ assigned the same label $(c_1,c_2,c_3)$. Two of $c_1,c_2,c_3$ are the same. That is, for some distinct $i,j \in \{1,2,3\}$, the edge set $\{ v_iw, v_jw: w \in W \}$ induces a monochromatic $K_{2,\lceil n/8\rceil}$. 
\end{proof}

Note that several of our lemmas include the condition $N_{\torso{G}{B_x}}(\apex(\Phi_x))\subseteq\topp(\Phi_x)$ for an almost-embedding $\Phi_x=(G^x_0,G^x_1,\ldots,G^x_p,A^x)$ of a torso $\torso{G}{B_x}$.  Since $B_x\setminus\topp(\Phi_x)=V(G^x_0)\setminus V(G^x_1\cup\cdots\cup G^x_p)$ and $\torso{G}{B_x}-A^x=G^x_0\cup G^x_1\cup\cdots\cup G^x_p$, this condition is equivalent to $N_{\torso{G}{B_x}}(B_x\setminus\topp(\Phi_x))\subseteq V(G^x_0)$ (that is, $N_{\torso{G}{B_x}}(v)\subseteq V(G^x_0)$ for every $v\in B_x\setminus\topp(\Phi_x)$).

\begin{lem}\label{K3nBound}
  Let $g,p,k,a,\Delta\in\NN$, let $G$ be a graph, and let $(\calT,\calE):=((B_x:x\in V(T)), (\Phi_x:x\in V(T)))$ be a rooted $(g,p,k,a)$-almost-embedded tree-decomposition of $G$ such that
  \begin{enumerate}[label=\rm(\alph*)]
    \item for each $x\in V(T)$, $N_{\torso{G}{B_x}}(\apex(\Phi_x))\subseteq\topp(\Phi_x)$;
    \item for each node $y$ of $T$ with parent $x$, 
    \begin{enumerate}[label=\rm(\roman*)]
      \item $B_x\cap B_y\subseteq\apex(\Phi_y)$ and 
      \item $(B_x\cap B_y)\cap\topp(\Phi_x)\neq\emptyset$; and
    \end{enumerate}
    \item $\deg_G(v)\le\Delta$ for each $v\in \bigcup_{x\in V(T)}\topp(\Phi_x)$.
  \end{enumerate}
  Then $G$ has no $K_{3,n}$-subdivision, where $n:= \max\{\Delta+1,\mks+1,a+24(2g+1),a+8p(2k+5)(2g+1)\}$. 
\end{lem}

\begin{proof}
  Suppose for the sake of contradiction that $G$ contains a $K_{3,n}$-subdivision $Q$ with degree-$n$ vertices $v_1,v_2,v_3$ and degree-$3$ vertices $w_1,\ldots,w_n$ (note that $n\ge 24(2g+1)\ge 24\ge 4$).
  
  For each $x\in V(T)$, let $\Phi_x:=(G^x_0,G^x_1,\ldots,G^x_p,A^x)$ and let $\Sigma^x$ be the surface on which $G^x_0$ is embedded.
  Define $\topp(G):=\bigcup_{x\in V(T)}\topp(\Phi_x)$.  If $v_i\in\topp(G)$ for some $i\in\{1,2,3\}$, then $n=\deg_Q(v_i)\le\deg_G(v_i) \le \Delta$ by~(c), contradicting $n\ge\Delta+1$.  Thus $v_1,v_2,v_3\not\in\topp(G)$.  For each non-root node $y$ of $T$ with parent $x$, $B_x\cap B_y\subseteq\apex(\Phi_y)\subseteq\topp(G)$ by~(b)(i), so the adhesion $B_x\cap B_y$ is contained in $\topp(G)$.  Therefore, $v_1,v_2,v_3$ do not belong to any adhesion, which implies that, for each $i\in\{1,2,3\}$, there exists a unique $x_i\in V(T)$ such that $v_i\in B_{x_i}$.

  Suppose that $x_i\neq x_j$ for some $i,j\in\{1,2,3\}$.  Let $xy$ be an edge of $T$ on the path from $x_i$ to $x_j$.  Then $\{v_i,v_j\}\cap B_x\cap B_y= \emptyset$. 
  Then every path from $v_i$ to $v_j$ in $G$ contains an internal vertex in $B_x\cap B_y$.  Since $Q$ contains $n$ pairwise internally disjoint paths from $v_i$ to $v_j$, this implies that $n\le |B_x\cap B_y|\le \mks$ by \cref{mks}, contradicting $n\ge\mks+1$.  Thus $x:=x_1=x_2=x_3$, so $v_1,v_2,v_3\in B_x\setminus\topp(\Phi_x)\subseteq V(G^x_0)$ and $v_1,v_2,v_3\not\in B_y$ for all $y\in V(T)\setminus\{x\}$. 

  For each $i\in\{1,2,3\}$ and each $j\in\{1,\ldots,n\}$, let $Q_{i,j}$ be the path in $Q$ between $v_i$ and $w_j$. For each $j\in\{1,\ldots,n\}$, let $S_j:=\bigcup_{i=1}^3 (Q_{i,j}-v_i)$.  Since the graphs $S_1,\ldots,S_n$ are pairwise vertex-disjoint and $|A^x|\le a$, there are at most $a$ indices $j$ such that $S_j$ contains a vertex in $A^x$.  Without loss of generality, assume that $S_1,\ldots,S_{n-a}$ contain no vertices of $A^x$.  Note that $n-a\ge 24(2g+1)\ge 24\ge 1$.  For each $i\in\{1,2,3\}$ and each $j\in\{1,\ldots,n-a\}$, say that $Q_{i,j}$ is \emph{red} if $Q_{i,j}\subseteq G^x_0$ and \emph{blue} otherwise.  By \cref{MonoK2n}, $Q$ contains a monochromatic $K_{2,\lceil (n-a)/8\rceil}$-subdivision $Q'$.  Without loss of generality, $Q':=\bigcup_{i=1}^2\bigcup_{j=1}^{\lceil (n-a)/8\rceil}Q_{i,j}$.
  
  First suppose that the paths $Q_{i,j}$ are red for all $i\in\{1,2\}$ and all $j\in\{1,\ldots,\lceil (n-a)/8\rceil\}$.  That is, $Q'$ is a subgraph of the surface-embedded graph $G^x_0$.  Since $n\ge a + 24(2g+1)$, we have $\lceil (n-a)/8\rceil \geq 3(2g+1) \ge \max\{3, 2(2g+1)\}$.  By \cref{homotopic_disc} with $c=3$ applied to the $v_1v_2$-paths of $Q'$, $Q'$ contains a $K_{2,3}$-subdivision $Q''$ with degree-$3$ vertices $v_1$ and $v_2$ such that $Q''$ is contained in a closed disc $D\subset \Sigma^x$ whose boundary is the union of two of the three $v_1v_2$-paths of $Q''$. Refer to \cref{k3nbound_fig}(a).  Without loss of generality, $Q'':=\bigcup_{i=1}^2\bigcup_{j=1}^3 Q_{i,j}$.  The embedding of $Q''$ on $\Sigma^x$ has three faces, two of which are open discs contained in $D$ and the third of which is $\Sigma^x\setminus D$.  Each face of $Q''$ is bounded by the union of two of the three $v_1v_2$-paths of $Q''$, so its closure contains two vertices of $\{w_1,w_2,w_3\}$ and omits the third.  Let $f$ be the face of $Q''$ that contains $v_3$, which is well-defined since $v_3\in V(G^x_0)\setminus V(Q'')\subseteq\Sigma^x\setminus Q''$.  Without loss of generality, assume that $w_1$ and $w_2$ are on the boundary of $f$, so $w_3$ is not contained in the closure of $f$.  

  \begin{figure}
    \begin{center}
      \begin{tabular}{cc}
        \includegraphics[page=1]{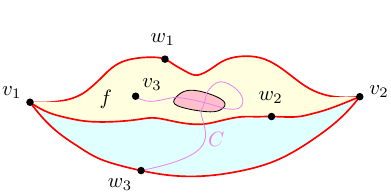} &
        \includegraphics[page=2]{figs/k3nbound} \\
        (a) & (b)
      \end{tabular}
      \caption{Two cases in the proof of \cref{K3nBound}}
      \label{k3nbound_fig}
    \end{center}
  \end{figure}
  
  Let $P_0:=Q_{3,3}$ be the path from $v_3$ to $w_3$ in $Q$; note that $v_3,w_3\in V(G^x_0)\subseteq B_x$. Consider a maximal subpath $P'$ of $P_0$ of length at least $2$ whose internal vertices do not belong to $B_x$.  Then the endpoints of $P'$ are distinct vertices in $B_x\cap B_y$ for some neighbour $y$ of $x$ in $T$.  Replacing each such subpath $P'$ of $P_0$ with the direct edge in $\torso{G}{B_x}$ between its two endpoints yields a path $P$ in $\torso{G}{B_x}$ from $v_3$ to $w_3$ with $V(P)=V(P_0)\cap B_x$.  Since $V(P) \setminus \{v_3\} \subseteq V(S_3)$, $S_3$ contains no vertices of $A^x$, and $v_3\notin A^x$, the path $P$ is disjoint from $A^x$ and satisfies $V(P)\cap V(Q'')=\{w_3\}$. Thus, $P$ is a path in $\torso{G}{B_x} - A^x = G^x_0 \cup G^x_1 \cup \dots \cup G^x_p$. Consider the curve $C\subseteq\Sigma^x$ obtained from $P$ as follows: $C$ contains each edge in $E(P)\cap E(G^x_0)$.  Since the graphs $G^x_1, \dots, G^x_p$ are pairwise vertex-disjoint and $v_3,w_3\in V(G^x_0)$, each non-trivial component $R$ of $P-E(G^x_0)$ is a path in $\torso{G}{B_x}$ whose edges are contained in $E(G^x_i)\setminus E(G^x_0)$ for a unique $i\in\{1,\ldots,p\}$ and whose endpoints are distinct vertices in $V(G^x_0)\cap V(G^x_i)\subseteq\partial D^x_i$. The curve $C$ replaces each such path $R$ with an arc in $D^x_i$ between the endpoints of $R$ whose interior is contained in $\operatorname{int}(D^x_i)$. Since each $D^x_i$ is $G^x_0$-clean and $V(P)\cap V(Q'')=\{w_3\}$, $C$ is a (possibly self-intersecting) curve in $\Sigma^x$ from $v_3$ to $w_3$ such that $C\setminus\{w_3\}$ is connected and disjoint from $Q''$.  Since $v_3\in f$, it follows that $C\setminus\{w_3\}\subseteq f$ and hence $w_3$ lies in the closure of $f$, a contradiction.
  
  Now assume that $Q'$ is blue.  That is, $Q_{i,j}$ is not a subgraph of $G^x_0$ for each $i\in\{1,2\}$ and each $j\in\{1,\ldots,\lceil (n-a)/8\rceil\}$.  For each $j\in\{1,\ldots,\lceil (n-a)/8\rceil\}$, $V(Q_{1,j})\not\subseteq B_x\setminus\topp(\Phi_x)$ since, otherwise, $Q_{1,j}\subseteq G[B_x\setminus\topp(\Phi_x)]\subseteq G^x_0$.  Since $v_1\in B_x\setminus\topp(\Phi_x)$, let $z_j$ be the first vertex along $Q_{1,j}$ (from $v_1$ to $w_j$) that does not belong to $B_x\setminus\topp(\Phi_x)$.  Then $z_j\neq v_1$, and every vertex of $Q_{1,j}$ strictly before $z_j$ belongs to $B_x\setminus\topp(\Phi_x)\subseteq V(G^x_0)$.
  
  If $z_j\in B_x$, then since $z_j\in V(S_j)$ and $S_j\cap A^x=\emptyset$, we have $z_j\in\topp(\Phi_x)\setminus A^x=V(G^x_1\cup\cdots\cup G^x_p)$.  Moreover, the edge of $Q_{1,j}$ entering $z_j$ has its other endpoint in $B_x\setminus\topp(\Phi_x)$ and therefore belongs to $E(G^x_0)$.  Thus $z_j\in V(G^x_0)\cap V(G^x_{\ell_j})$ for some $\ell_j\in\{1,\ldots,p\}$; in this case, we set $a_j:=z_j$ and $a'_j:=z_j$.  Otherwise, $z_j\not\in B_x$, and we let $a_j$ be the predecessor of $z_j$ on $Q_{1,j}$.  Since $a_jz_j\in E(G)$ with $a_j\in B_x$ and $z_j\not\in B_x$, we have $a_j\in B_x\cap B_y$ for some neighbour $y$ of $x$ in $T$; in particular, $a_j\neq v_1$ because $v_1$ belongs to no adhesion of $\calT$.  Since $a_j\not\in\topp(\Phi_x)$, $y$ cannot be the parent of $x$ (since condition (b)(i) would imply that $a_j\in\apex(\Phi_x)\subseteq \topp(\Phi_x)$), so $y$ is a child of $x$.  By~(b)(ii), there exists a vertex $u\in (B_x\cap B_y)\cap\topp(\Phi_x)$. Then $u\in N_{\torso{G}{B_x}}(a_j)$ since $B_x\cap B_y$ is a clique in $\torso{G}{B_x}$ and $u\neq a_j$.  Since $N_{\torso{G}{B_x}}(\apex(\Phi_x))\subseteq\topp(\Phi_x)$ by~(a) and $a_j\in B_x\setminus\topp(\Phi_x)$, we have $N_{\torso{G}{B_x}}(a_j)\subseteq V(G^x_0)$.  Thus $u\in V(G^x_0)\cap (\topp(\Phi_x)\setminus A^x)=V(G^x_0)\cap V(G^x_1\cup\cdots\cup G^x_p)$, so $u\in V(G^x_0)\cap V(G^x_{\ell_j})$ for some $\ell_j\in\{1,\ldots,p\}$.  Since $a_ju\in E(\torso{G}{B_x}-A^x)$ and $a_j\not\in\topp(\Phi_x)$, the edge $a_ju$ belongs to $E(G^x_0)$. In this case, we set $a'_j:=u$.  
  
  Since $\lceil (n-a)/8\rceil \ge 1$, $p\ge 1$ (the blue case cannot occur when $p=0$).  Since $a_j\in V(S_j)$ for each $j$, the vertices $a_1,\ldots,a_{\lceil (n-a)/8\rceil}$ are pairwise distinct, though $a'_1,\ldots,a'_{\lceil(n-a)/8\rceil}$ are not necessarily distinct.  Since $a'_j\in V(G^x_0)\cap V(G^x_{\ell_j})$ for each $j\in\{1,\ldots,\lceil(n-a)/8\rceil\}$, by the Pigeonhole Principle we may assume without loss of generality that $a'_j\in V(G^x_0)\cap V(G^x_1)$ for each $j\in\{1,\ldots,q\}$, where $q:=\lceil (n-a)/8p\rceil$.  For each $j\in\{1,\ldots,q\}$, let $P_j$ be the subpath of $Q_{1,j}$ from $v_1$ to $a_j$, and let $P'_j:=P_j$ if $a'_j=a_j$, or let $P'_j$ be obtained by appending the edge $a_ja'_j$ to $P_j$ if $a'_j\neq a_j$.  By construction, each $P'_j$ is a path in $G^x_0$ from $v_1$ to $a'_j$ with $(V(P'_j)\setminus\{a'_j\})\cap\topp(\Phi_x)=\emptyset$.  Consequently, for distinct $j,j'\in\{1,\ldots,q\}$, the paths $P'_j$ and $P'_{j'}$ are edge-disjoint 
  and intersect only at $v_1$ and possibly at a shared endpoint $a'_j=a'_{j'}$.  Let $D:=D^x_1$ be the $G^x_0$-clean disc whose boundary contains $V(G^x_0)\cap V(G^x_1)$, including $U:=\{a'_1,\ldots,a'_q\}$.  
  
  Refer to \cref{k3nbound_fig}(b).
  For each $u\in U$, choose a closed disc $N_u\subseteq\Sigma^x$ containing $u$ in its interior such that the discs $(N_u:u\in U)$ are pairwise disjoint, each $N_u$ is disjoint from $D^x_2,\ldots,D^x_p$, from all vertices of $G^x_0\setminus\{u\}$, and from all edges of $G^x_0$ not incident to $u$, each edge of $G^x_0$ incident to $u$ intersects $N_u$ in a single radial arc from $\partial N_u\setminus D$ to $u$, and $N_u\cap D$ is a closed disc with $N_u\cap\partial D$ a single arc containing $u$ in its interior.  Place a point $z$ in $\operatorname{int}(D)\setminus\bigcup_{u\in U} N_u$.  For each $u\in U$ and those $j\in\{1,\ldots,q\}$ with $a'_j=u$, replace $P'_j\cap N_u$ by pairwise disjoint arcs from $P'_j\cap\partial N_u$ to distinct points $q_j\in\operatorname{int}(N_u)\cap\partial D$ whose interiors lie in $\operatorname{int}(N_u)\setminus D$ (preserving their order around $\partial N_u$), and extend them from $q_j$ to $z$ by pairwise internally disjoint arcs in $\operatorname{int}(D)$ avoiding $\bigcup_{u'\in U\setminus\{u\}} N_{u'}$ so that the resulting curves $P''_1,\ldots,P''_q$ from $v_1$ to $z$ have pairwise disjoint interiors.  Let $X:=\bigcup_{j=1}^{q}P''_j$, so $X$ is a $K_{2,q}$-subdivision embedded in $\Sigma^x$.
  
  Let $r:=2k+5\ge 5$.  Since $n\ge a+8p(2k+5)(2g+1)$, we have $q=\lceil (n-a)/8p\rceil \ge r(2g+1) \ge \max\{r, (r-1)(2g+1)\}$.  By \cref{homotopic_disc} with $c=r$, $X$ contains a $K_{2,r}$-subdivision $X'$ that is contained in a closed disc $F\subset \Sigma^x$ whose boundary is the union of two of the $v_1z$-paths in $X'$.  Without loss of generality, assume that $X'=\bigcup_{j=1}^r P''_j$.  Then $X'$ has $r$ faces in $\Sigma^x$, $r-1$ of which are open discs contained in $F$ and one of which is $\Sigma^x\setminus F$.  Each face of $X'$ is bounded by the union of two curves in $\{P''_1,\ldots,P''_r\}$.  Since $v_2\in V(G^x_0)\setminus V(P'_1\cup\cdots\cup P'_r)$ and $v_2\notin\topp(\Phi_x)$, we have $v_2\in\Sigma^x\setminus X'$; let $f$ be the face of $X'$ containing $v_2$.  Without loss of generality, assume that $\partial f = P''_1\cup P''_2$, so $P''_3\setminus\{v_1,z\},\ldots,P''_r\setminus\{v_1,z\}$ are disjoint from the closure $\overline{f}$ of $f$.  Since $a'_1,a'_2\in V(G^x_0)\cap V(G^x_1)$, the path-decomposition of $G^x_1$ has bags $B_1$ and $B_2$ indexed by $a'_1$ and $a'_2$, respectively (with $B_1=B_2$ if $a'_1=a'_2$). 
  
  We claim that $Z:=B_1\cup B_2\cup V(P_1\cup P_2)$ separates $v_2$ from $\{a_3,\ldots,a_r\}$ in $\torso{G}{B_x} - A^x$.  Note that $V(P'_1\cup P'_2)\subseteq Z$ since $a'_1\in B_1$ and $a'_2\in B_2$.  Say that a vertex $v$ of $\torso{G}{B_x}-(A^x \cup Z)$ is an \emph{inner vertex} if either $v\in V(G^x_0)\cap f$, or $v\in V(G^x_i)\setminus V(G^x_0)$ for some $i\in\{1,\ldots,p\}$ and $v$ belongs to a bag of the path-decomposition of $G^x_i$ indexed by a vertex in $V(G^x_0)\cap V(G^x_i)\cap f$.  Clearly $v_2$ is an inner vertex.  For each $j\in\{3,\ldots,r\}$, if $a_j\notin Z$, then $a_j\notin\{a'_1,a'_2\}$; either $a_j=a'_j\in U\setminus\{a'_1,a'_2\}$ (in which case $N_{a'_j}$ is disjoint from $P''_1\cup P''_2=\partial f$ and meets $P''_j\setminus\{v_1,z\}\subseteq\Sigma^x\setminus\overline{f}$, so $a_j\in N_{a'_j}\subseteq\Sigma^x\setminus\overline{f}$), or $a_j\neq a'_j$ (in which case $a_j\in P'_j\setminus (N_{a'_j}\cup D)\subseteq P''_j\setminus\{v_1,z\}\subseteq\Sigma^x\setminus\overline{f}$).  Thus none of $a_3,\ldots,a_r$ is an inner vertex.

  We next verify that for every inner vertex $v$, $N_{\torso{G}{B_x}}(v)$ is contained in the union of $A^x\cup Z$ and the set of inner vertices.  First, if $vw\in E(G^x_0)$ with $v\in V(G^x_0)\cap f$ and $w\notin Z$, then $v,w\notin V(P'_1\cup P'_2)$ (in particular $v,w\notin\{a'_1,a'_2\}$), so the embedded edge $vw$ is disjoint from $P'_1\cup P'_2$, from $N_{a'_1}\cup N_{a'_2}$, and from $\operatorname{int}(D)$, and hence from $P''_1\cup P''_2=\partial f$; thus $w\in f$ is an inner vertex.  Second, along $\partial D$, the two distinct points $q_1,q_2$ where $P''_1$ and $P''_2$ cross $\partial D$ partition $(V(G^x_0)\cap V(G^x_1))\setminus\{a'_1,a'_2\}$ into the subset lying in $f$ and the subset lying in $\Sigma^x\setminus\overline{f}$ (since $P''_3\cap\operatorname{int}(D)$ lies outside $\overline{f}$ while $z\in\partial f$).  Because no vertex of $V(G^x_0)\setminus\{a'_1,a'_2\}$ lies in $N_{a'_1}\cup N_{a'_2}$, removing $\{a'_1,a'_2\}$ from $\partial D$ induces the exact same partition of $(V(G^x_0)\cap V(G^x_1))\setminus\{a'_1,a'_2\}$ (with one part empty if $a'_1=a'_2$).  Consequently, in the path-decomposition of $G^x_1$, every path of bags between a bag indexed by a vertex in $f$ and a bag indexed by a vertex in $\Sigma^x\setminus\overline{f}$ passes through $B_1$ or $B_2$; by running intersection, no vertex of $V(G^x_1)\setminus (B_1\cup B_2)$ belongs to bags of both types, so every edge of $G^x_1-(B_1\cup B_2)$ incident to an inner vertex has its other endpoint inner as well.  Third, for each $i\in\{2,\ldots,p\}$, since $(V(P'_1\cup P'_2)\setminus\{a'_1,a'_2\})\cap\topp(\Phi_x)=\emptyset$, the $G^x_0$-clean disc $D^x_i$ is disjoint from $P''_1\cup P''_2=\partial f$, so either $D^x_i\subseteq f$ (and every vertex of $G^x_i$ is inner) or $D^x_i\subseteq\Sigma^x\setminus\overline{f}$ (and no vertex of $G^x_i$ is inner).  Hence the component of $\torso{G}{B_x}-(A^x\cup Z)$ containing $v_2$ consists solely of inner vertices and contains none of $a_3,\ldots,a_r$, proving that $Z$ separates $v_2$ from $\{a_3,\ldots,a_r\}$ in $\torso{G}{B_x}-A^x$.
  
  For each $j\in\{3,\ldots,r\}$, let $R_j:=Q_{1,j}\cup Q_{2,j}$ be the path in $Q$ from $v_1$ to $v_2$ that includes $w_j$, and let $R^*_j$ be the subpath of $R_j$ from $a_j$ to $v_2$.  Since $V(R^*_j)\subseteq V(S_j)\cup\{v_2\}$ is disjoint from $A^x$ and $a_j,v_2\in B_x$, replacing maximal subpaths of $R^*_j$ of length at least $2$ whose internal vertices lie outside $B_x$ by torso edges yields a path in $\torso{G}{B_x}-A^x$ from $a_j$ to $v_2$ whose vertex set is $V(R^*_j)\cap B_x$.  Thus $R^*_j$ contains a vertex $u_j\in Z$.  Since $v_1\notin V(R^*_j)$ (as $a_j\neq v_1$) and $v_2\notin Z$, the vertex $u_j$ is an internal vertex of $R_j$.  Since $R_1,\ldots,R_r$ are pairwise internally disjoint and $V(P_1\cup P_2)\subseteq V(R_1\cup R_2)$, we have $u_j\in (B_1\cup B_2)\setminus V(P_1\cup P_2)\subseteq B_1\cup B_2$ for each $j\in\{3,\ldots,r\}$, and $u_3,\ldots,u_r$ are pairwise distinct.  Therefore, $2k+3 = r-2\le |B_1\cup B_2| \le 2(k+1)$, a contradiction.
\end{proof}

\section{The Proof}
\label{proof}

In this section, we complete the proof of \cref{main}. As outlined in the introduction, the proof partitions the tree $T$ used in an almost-embedded tree-decomposition into subtrees.  Within each subtree, each adhesion $B_x\cap B_y$ between a node $x$ and its child $y$ is contained in the apex set of $\torso{G}{B_y}$ and intersects the top of $\torso{G}{B_x}$. The following lemma shows how to handle one of these subtrees:

\begin{lem}\label{base_case}
  There exists a function $f:\NN^5\to\NN^+$ such that the following is true, for all $g,p,k,a,\Delta\in\NN$:  Let $G$ be a graph of maximum degree at most $\Delta$ and let $(\calT,\calE):=((B_x:x\in V(T)),(\Phi_x:x\in V(T)))$ be a rooted $(g,p,k,a)$-almost-embedded tree-decomposition of $G$ such that 
  \begin{enumerate}[label=\rm(\alph*)]
    \item for each $x\in V(T)$, $N_{\torso{G}{B_x}}(\apex(\Phi_x))\subseteq \topp(\Phi_x)$;
    \item for each node $y$ of $T$ with parent $x$,
    \begin{enumerate}[label=\rm(\roman*)]
      \item\label{apex_adhesions} $B_x\cap B_y\subseteq \apex(\Phi_y)$, and
      \item\label{top_adhesions_nonempty} $(B_x\cap B_y)\cap\topp(\Phi_x)\neq\emptyset$, 
    \end{enumerate}
  \end{enumerate}
  and let $S\subseteq \apex(\Phi_{r})$, where $r$ is the root of $T$.
  
  Then there exists a graph $H$ of treewidth at most $3$ and a $2$-weak layered $H$-partition $(\calL:=(L_1,L_2,\ldots),\calP:=(E_x:x\in V(H)))$ of $G$ of layered width at most $f(g,p,k,a,\Delta)$, such that $S=E_{x_0}\subseteq L_1$ for some $x_0\in V(H)$ with $\deg_H(x_0)=1$.
\end{lem}

\begin{proof}
  For each $x\in V(T)$, let $\Phi_x:=(G^x_0,G^x_1,\ldots,G^x_p,A^x)$ and let $\Sigma^x$ be the surface of Euler genus at most $g$ on which $G^x_0$ is embedded. For each $i\in\{1,\ldots,p\}$, let $D^x_i$ be the $G^x_0$-clean disc whose boundary contains the vortex boundary $V(G^x_i)\cap V(G^x_0)$.  Let $N:=N_{G}(S)$, and observe that $|N|\le\Delta |S|\le\Delta|A^{r}|\le\Delta a$.  We claim that $(S\cup N)\cap B_x\subseteq\topp(\Phi_x)$ for every $x\in V(T)$.  First, for each $s\in S\cap B_x$, if $x=r$ then $s\in\apex(\Phi_r)\subseteq\topp(\Phi_r)$, and if $x\neq r$ with parent $p$, then since $s\in B_r\cap B_x$, \ref{apex_adhesions} implies that $s\in B_p\cap B_x\subseteq\apex(\Phi_x)\subseteq\topp(\Phi_x)$; thus $S\cap B_x\subseteq\apex(\Phi_x)\subseteq\topp(\Phi_x)$.  Next, for each $w\in N\cap B_x$, choose $s\in S$ and $z\in V(T)$ such that $sw\in E(G[B_z])$ (note that $s\neq w$ since $S\cap N=\emptyset$).  If $z=x$, then $s\in S\cap B_x\subseteq\apex(\Phi_x)$ and $sw\in E(\torso{G}{B_x})$, so $w\in N_{\torso{G}{B_x}}[\apex(\Phi_x)]\subseteq\topp(\Phi_x)$ by~(a).  If $z\notin V(T_x)$, then $x\neq r$ and $w\in B_x\cap B_z$ belongs to the adhesion between $x$ and its parent, so $w\in\apex(\Phi_x)\subseteq\topp(\Phi_x)$ by \ref{apex_adhesions}.  Finally, if $z\in V(T_x-x)$, let $c$ be the child of $x$ on the path from $x$ to $z$ in $T$; since $s\in B_r\cap B_z$ and $w\in B_x\cap B_z$, $\{s,w\}\subseteq B_x\cap B_c$, so $sw\in E(\torso{G}{B_x})$ with $s\in\apex(\Phi_x)$ and hence $w\in N_{\torso{G}{B_x}}[\apex(\Phi_x)]\subseteq\topp(\Phi_x)$ by~(a).
  
  By \cref{almost_embedded_layered_partition}, for each $x\in V(T)$, there exists a graph $I^x$ of treewidth at most $\ell$ such that $\torso{G}{B_x}$ has a minimal connected layered $I^x$-partition $(\calL^x,\calP^x):=((L^x_1,L^x_2,\ldots),(C^x_y:y\in V(I^x)))$ of layered width at most $\ell$, where $\ell:=\ell(g,p,k,a)$, with $\topp(\Phi_x)\subseteq L^x_1$ and such that $C^x_v=\{v\}$ for each $v\in L^x_1\cup L^x_2$.  Choose the node sets indexing non-singleton parts of $\calP^x$ so that the sets $V(I^x)\setminus B_x$ ($x\in V(T)$) are pairwise disjoint and disjoint from $V(G)$.  See \cref{curtain_fig}.  

  \begin{figure}
    \begin{center}
      \begin{tabular}{ccc}
        \includegraphics[page=1]{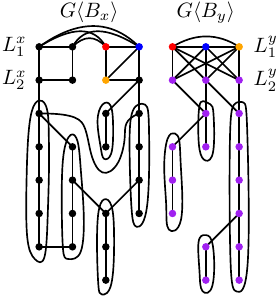} &
        \includegraphics[page=2,trim={10 0 10 0},clip]{figs/curtain} &
        \includegraphics[page=3,trim={10 0 10 0},clip]{figs/curtain}
      \end{tabular}
    \end{center}
    \caption{The torsos $\torso{G}{B_x}$ and $\torso{G}{B_y}$ (left), the quotient graphs $I^x$ and $I^y$ (middle), and $I=I^x\cup I^y$ (right) for a $2$-node tree $T$ where $x$ is the parent of $y$.}
    \label{curtain_fig}
  \end{figure}

  We will argue that the partitions $\calP^x$, for $x\in V(T)$, are compatible in the following sense:  Each vertex $v$ of $G$ that appears in an adhesion of $\calT$ is a singleton $C^x_v=\{v\}$ in each $\calP^x$ for which $v\in B_x$.  Consider an adhesion $B_x\cap B_y$ defined by some edge $xy$ of $T$ where $x$ is the parent of $y$, and let $v\in B_x\cap B_y$.  By \ref{apex_adhesions}, $v\in\apex(\Phi_y)\subseteq\topp(\Phi_y)\subseteq L^y_1$.  By \ref{top_adhesions_nonempty}, $B_x\cap B_y$ contains a vertex $u\in\topp(\Phi_x)\subseteq L^x_1$. Since $B_x\cap B_y$ is a clique in $\torso{G}{B_x}$ and $\calL^x$ is a layering of $\torso{G}{B_x}$, $v\in N_{\torso{G}{B_x}}[u]\subseteq L^x_1\cup L^x_2$.  Thus, for each adhesion $B_x\cap B_y$ of $\calT$ and each $v\in B_x\cap B_y$, $C^y_v=C^x_v=\{v\}\subseteq (L^x_1\cup L^x_2)\cap (L^y_1\cup L^y_2)$.  Consequently, $V(I^x)\cap V(I^z)=B_x\cap B_z$ for all distinct $x,z\in V(T)$.

  Let $I:=\bigcup_{x\in V(T)}I^x$. Then $\calI:=(V(I^x):x\in V(T))$ is a rooted tree-decomposition of $I$. For each edge $xy$ of $T$, $B_x\cap B_y=V(I^x)\cap V(I^y)$ is a clique in $\torso{G}{B_x}$.  Since $\calP^x$ is an $I^x$-partition of $\torso{G}{B_x}$ in which each vertex $v\in B_x\cap B_y$ forms a singleton part $C^x_v=\{v\}$, $V(I^x)\cap V(I^y)$ is also a clique in $I^x$.  Therefore $\torso{I}{V(I^x)}=I^x$ for each $x\in V(T)$.  Thus $\calI$ is a tree-decomposition of $I$ in which each torso has treewidth at most $\ell$, so $I$ has treewidth at most $\ell$.

  We now construct an $I$-partition of $G$.\footnote{In fact we construct something stronger: an $I$-partition of $\bigcup_{x\in V(T)}\torso{G}{B_x}\supseteq G$.}  For each $q\in V(I)$, choose a node $x\in V(T)$ with $q\in V(I^x)$ and define $C_q:=C^x_q$.  (Note that the choice of $x$ is irrelevant; if $q\in V(I^x)\cap V(I^y)$ for some $x\neq y$, then $q\in V(G)$ and $C^x_q=C^y_q=\{q\}$.)   Let $\calR:=(C_q:q\in V(I))$. We claim that $\calR$ is an $I$-partition of $G$.  First we establish that $\calR$ covers $V(G)$ with pairwise-disjoint sets.  Clearly, $V(G)=\bigcup_{x\in V(T)}B_x=\bigcup_{x\in V(T)}\bigcup_{q\in V(I^x)}C^x_q$, so $\calR$ covers $V(G)$.  To see why the sets in $\calR$ are disjoint, suppose that $v\in C^x_q\cap C^y_r$.  If $x=y$, then $q=r$ since $\calP^x$ is an $I^x$-partition of $\torso{G}{B_x}$.  If $x\neq y$, then $v\in B_x\cap B_y$, so $C^x_q=C^y_r=\{v\}=C_q$ (and $q=r=v$). Thus, for each $v\in V(G)$ there is a unique $q\in V(I)$ such that $C_q$ contains $v$.  Finally, we show that $I$ contains the quotient graph $G/\calR$.  For each edge $vw$ of $G$, $v,w\in B_x$ for some $x\in V(T)$. Therefore $vw$ is an edge of $G[B_x]\subseteq \torso{G}{B_x}$.  Then $v\in C^x_q$ and $w\in C^x_r$ for some $q,r\in V(I^x)$.
  Thus, $q=r$ or $qr\in E(I^x)$ since $\calP^x$ is an $I^x$-partition of $\torso{G}{B_x}$. 

  Thus far, we have constructed an $I$-partition $\calR$ of $G$ and a tree-decomposition $\calI$ of $I$ whose torsos have treewidth at most $\ell$.  We now show that each torso of $\calI$ has a $(g,p,k,a)$-almost-embedding. 
  Since $(\calL^x,\calP^x)$ is minimal and connected, and every non-singleton part of $\calP^x$ is disjoint from $\topp(\Phi_x)$ (and every edge of $\torso{G}{B_x}$ incident to $B_x\setminus\topp(\Phi_x)$ belongs to $E(G^x_0)$ by~(a)), the graph $I^x$ is obtained from $\torso{G}{B_x}$ by contracting connected subgraphs of $G^x_0-\topp(\Phi_x)$. These contractions can be done in such a way that $I^x$ has a $(g,p,k,a)$-almost-embedding $\Psi_x:=(I^x_0,I^x_1,\ldots,I^x_p,A^x)$ with $I^x_i=G^x_i$ for each $i\in\{1,\ldots,p\}$, in which $I^x_0$ is embedded on the same surface $\Sigma^x$ on which $G^x_0$ is embedded.  Since none of these contractions includes vertices of $G^x_0$ that participate in vortices, these contractions can be done so that, for each $i\in\{1,\ldots,p\}$, the $G^x_0$-clean disc $D^x_i$ used for the vortex $G^x_i$ is also the $I^x_0$-clean disc used for the vortex $I^x_i$.

  The pair $(\calI,\Psi:=(\Psi_x:x\in V(T)))$ is a rooted $(g,p,k,a)$-almost-embedded tree-decomposition of $I$.  Our goal is to apply \cref{CampbellProduct} to a graph derived from $I$.  However, the graph $I$ does not satisfy the degree-bound conditions of \cref{K3nBound}, so $I$ may contain $K_{3,n}$-subdivisions for arbitrarily large values of $n$, which makes \cref{CampbellProduct} inapplicable.  Specifically, the edges of $\torso{G}{B_x}\cap I^x$ not present in $G[B_x]$, for $x\in V(T)$, lead to vertices in $\topp(\Psi_x)=\topp(\Phi_x)$ that may have unbounded degree in $I$. This justifies the graphs $\widehat{I}^x$ for $x\in V(T)$ and $\widehat{I}=\bigcup_{x\in V(T)}\widehat{I}^x$ that we define next.

  For each $x\in V(T)$, let $\widehat{I}^x$ be the spanning subgraph of $I^x$ that contains an edge $qr\in E(I^x)$ if and only if there is an edge $vw$ in $G$ with $v\in C_q$ and $w\in C_r$. (Alternatively, $\widehat{I}^x$ is the quotient graph $G[B_x]/\calP^x$.)
  Thus, $\calP^x$ is a minimal $\widehat{I}^x$-partition of $G[B_x]$.  
  For each edge $vw\in E(\torso{G}{B_x})\setminus E(G[B_x])$, the vertices $v$ and $w$ are contained in an adhesion $B_x\cap B_y$ of $\calT$. Thus, $v,w\in L^x_1\cup L^x_2$, $C^x_v=\{v\}$ and $C^x_w=\{w\}$ are singletons in $\calP^x$, and every edge of $E(I^x)\setminus E(\widehat{I}^x)$ has both endpoints in $B_x\cap B_y$ for some $xy\in E(T)$. Thus, for each part $C^x_q$ of $\calP^x$, $G[B_x][C^x_q]=\torso{G}{B_x}[C^x_q]$ is connected. Thus, $(\calL^x,\calP^x)$ is a minimal connected layered $\widehat{I}^x$-partition of $G[B_x]$.  Let $\widehat{I}:=\bigcup_{x\in V(T)}\widehat{I}^x$.  Then $\calR=(C_q:q\in V(\widehat{I}))$ is a minimal connected $\widehat{I}$-partition of $G$.  For each $x\in V(T)$, and each $v\in \topp(\Psi_x)$, $C_v=\{v\}$ is a singleton in $\calR$. By \cref{singleton_degree}, $\deg_{\widehat{I}}(v)\le \deg_G(v)$.

  \begin{figure}
    \begin{center}
      \begin{tabular}{ccc}
        \includegraphics[page=1]{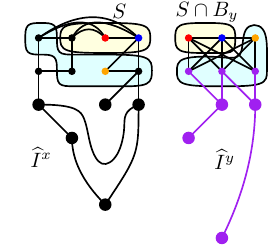} &
        \includegraphics[page=2]{figs/curtain2} &
        \includegraphics[page=3]{figs/curtain2}
      \end{tabular}
    \end{center}
    \caption{Constructing $J$ from $\widehat{I}^x$ and $\widehat{I}^y$ for a $2$-node tree $T$ where $x$ is the parent of $y$.}
    \label{curtain_fig2}
  \end{figure}

  We will now use \cref{K3nBound} to show that $\widehat{I}$ has no $K_{3,n}$-subdivision for a suitably large value of $n$. 
  Since $\widehat{I}$ is a spanning subgraph of $I$, the rooted tree-decomposition $\calI:=(V(I^x):x\in V(T))$ of $I$ is also a rooted tree-decomposition of $\widehat{I}$. For each $x\in V(T)$, since the adhesions of $\calI$ incident to $x$ are $V(I^x)\cap V(I^y)=B_x\cap B_y$ ($xy\in E(T)$), the definition of $\widehat{I}$ implies that $\torso{\widehat{I}}{V(I^x)}=\torso{I}{V(I^x)}=I^x$.
  Therefore, the pair $(\calI,\Psi:=(\Psi_x:x\in V(T)))$ is a rooted $(g,p,k,a)$-almost-embedded tree-decomposition of $\widehat{I}$. For each $x\in V(T)$, $\topp(\Psi_x)=\topp(\Phi_x)$, $\apex(\Psi_x)=\apex(\Phi_x)$, and $I^x$ is obtained from $\torso{G}{B_x}$ by contracting subgraphs in $B_x\setminus\topp(\Phi_x)$.  Each subgraph contracted this way has no vertex in $N_{\torso{G}{B_x}}(\apex(\Phi_x))$ by~(a), so assumptions~(a) and~(b) of the current lemma imply that $\widehat{I}$ and $(\calI,\Psi)$ satisfy requirements~(a) and~(b) of \cref{K3nBound}. As argued above, $\deg_{\widehat{I}}(v)\le\deg_{G}(v)\le\Delta$ for each $v\in \bigcup_{x\in V(T)}\topp(\Psi_x)$, so $\widehat{I}$ and $(\calI,\Psi)$ satisfy condition~(c) of \cref{K3nBound}. Therefore, \cref{K3nBound} implies that $\widehat{I}$ has no $K_{3,n}$-subdivision for $n:=\max\{\Delta+1, \mks+1, a+24(2g+1), a+8p(2k+5)(2g+1)\}$.
  
  Refer to \cref{curtain_fig2}.
  Introduce two distinct new vertices $y_0$ and $y_1$ not in $V(\widehat{I})$.  Fix some $x\in V(T)$. We modify $\calP^x$ and $\widehat{I}^x$ to obtain a graph $J^x$ and a $J^x$-partition $\calQ^x:=(C^x_y:y\in V(J^x))$ of $G[B_x]$ in which $C^x_{y_0}=S\cap B_x$ and $y_0$ is a vertex of degree $1$ in $J^x$ whose neighbour $y_1$ has $C^x_{y_1}=N\cap B_x$.  We do this in two steps.  Start with $J^x:=\widehat{I}^x$ and $\calQ^x:=\calP^x$.  In the first step, we remove $(S\cup N)\cap B_x$ from $J^x$. Since $(S\cup N)\cap B_x\subseteq\topp(\Phi_x)\subseteq L^x_1$, each vertex in $(S\cup N)\cap B_x$ forms a singleton part in $\calP^x$. Removing $(S\cup N)\cap B_x$ from $J^x$ therefore removes the parts $\{C^x_v=\{v\}:v\in (S\cup N)\cap B_x\}$ from $\calQ^x$, so $(\calL^x,\calQ^x)$ is a connected layered $J^x$-partition of $G[B_x]-(S\cup N)$.  In the second step, we add $y_0$ and $y_1$ to $J^x$ with $C^x_{y_0}:=S\cap B_x$ and $C^x_{y_1}:=N\cap B_x$, add the edge $y_0y_1$ to $J^x$, and make $y_1$ adjacent to every node $y\in V(J^x)\setminus\{y_0,y_1\}$ such that $G[B_x]$ contains an edge $vw$ with $v\in C^x_{y_1}$ and $w\in C^x_y$. Since $S\cap N=\emptyset$ and $N_G(S)=N$ (so no edge of $G$ joins $C^x_{y_0}$ to $B_x\setminus (S\cup N)$), $(\calL^x,\calQ^x)$ is a layered $J^x$-partition of $G[B_x]$.

  We emphasize that the vertices $y_0$ and $y_1$ are not local to $J^x$; these vertices appear in $J^{z}$ for every $z\in V(T)$. Let $C_{y_0}:=S$, let $C_{y_1}:=N$, and let $J:=\bigcup_{x\in V(T)} J^x$; note that $N_J(y_0)=\{y_1\}$.  
  Recall that, for each $x\in V(T)$ and each $y\in V(J^x)\setminus\{y_0,y_1\}$, $C_y=C^x_y$.    
  Since $\calQ^x$ is a $J^x$-partition of $G[B_x]$ for each $x\in V(T)$, $\calQ:=(C_y:y\in V(J))$ is a $J$-partition of $G$.
  Let $L'_i:=\bigcup_{x\in V(T)} L^x_i$ for each $i\in\NN^+$. For each $v\in V(G)$, let $\ell(v) := \min \{i\in\NN^+ : v \in L'_i\}$. Let $L_i := \{v \in V(G) : \ell(v) = i\}$ for each $i\in\NN^+$. Then $\calL:=(L_1,L_2,\ldots)$ is a sequence of pairwise disjoint subsets of $V(G)$ that cover $V(G)$. 
  
  We claim that $\calL$ is a $2$-weak layering of $G$. For each $v\in V(G)$, if there is a unique $x\in V(T)$ such that $v\in B_x$, then $v\in L^x_i$ implies $\ell(v)=i$. Otherwise, let $x$ be the node of $T$ closest to the root of $T$ such that $v\in B_x$.  Then $v\in B_x\cap B_y$ for some child $y$ of $x$, and $B_x\cap B_y$ contains a vertex $u\in\topp(\Phi_x)$ by \ref{top_adhesions_nonempty}. Therefore, $u\in L^x_1$ and $v\in N_{\torso{G}{B_x}}[u]\subseteq L^x_1\cup L^x_2$ since $B_x\cap B_y$ is a clique in $\torso{G}{B_x}$ and $\calL^x$ is a layering of $\torso{G}{B_x}$.  For each $y\in V(T)\setminus\{x\}$ with $v\in B_y$, $v$ belongs to the adhesion between $y$ and its parent, so $v\in\apex(\Phi_y)\subseteq\topp(\Phi_y)\subseteq L^y_1$, which implies $\ell(v)=1$.  Thus, for any $v\in V(G)$ and any $z\in V(T)$, if $v\in L^z_i$, then $\ell(v)\in\{i-1,i\}$.  For each $vw\in E(G)$, there exists $z\in V(T)$ with $v,w\in B_z$.  Since $\calL^z$ is a layering of $\torso{G}{B_z}$, we have $v\in L^z_i$ and $w\in L^z_j$ for some $i,j\in\NN^+$ with $|i-j|\le 1$. Then $\ell(v)\in\{i-1,i\}$ and $\ell(w)\in\{j-1,j\}$, so $|\ell(v)-\ell(w)|\le |i-j|+1\le 2$, as required by the definition of $2$-weak layering. 

  We now bound the layered width of $(\calL,\calQ)$. For each $y\in V(J)\setminus\{y_0,y_1\}$, $C_y=C^x_y$ for some $x\in V(T)$. For each $i\in\NN^+$, each vertex $v\in C_y\cap L_i$ satisfies $\ell(v)=i$. Since $v\in B_x$, $v\in L^x_j$ for some $j\in\NN^+$, and since $\ell(v)\in\{j-1,j\}$, we have $j\in\{i,i+1\}$. Therefore $C_y\cap L_i\subseteq (C^x_y\cap L^x_i)\cup (C^x_y\cap L^x_{i+1})$, so $|L_i\cap C_y|\le |C^x_y\cap L^x_i| + |C^x_y\cap L^x_{i+1}|\le 2\ell$. Since $S\subseteq\apex(\Phi_r)\subseteq L^r_1$ and $(S\cup N)\cap B_x\subseteq\topp(\Phi_x)\subseteq L^x_1$ for each $x\in V(T)$, we have $C_{y_0}\cup C_{y_1}=S\cup N\subseteq L_1$, so $L_i\cap C_{y_0}=L_i\cap C_{y_1}=\emptyset$ for each $i\in\NN^+\setminus\{1\}$.  Furthermore, $|C_{y_0}\cap L_1|=|S|\le a$ and $|C_{y_1}\cap L_1|=|N|\le \Delta a$.  Thus, $(\calL,\calQ)$ is a $2$-weak layered $J$-partition of $G$ of layered width at most $\ell':=\max\{2\ell,a,\Delta a\}$.  

  Set $n':=1+\max\{\Delta^2 a,n\}\ge n+1\ge 25$; we claim that $J-y_0$ contains no $K_{3,n'}$-subdivision.  First consider the graph $J_2:=J-\{y_0,y_1\}$.  By definition, $J_2=\widehat{I}-(S\cup N)$, so $J_2$ is a subgraph of $\widehat{I}$, and hence $J_2$ has no $K_{3,n}$-subdivision.  Suppose for a contradiction that $J-y_0$ contains a $K_{3,n'}$-subdivision $Q$. If $y_1\notin V(Q)$, then $Q$ is a subgraph of $J_2$, contradicting $n'>n$. If $y_1\in V(Q)$ and $\deg_Q(y_1)=n'$, then $n'=\deg_Q(y_1)\le\deg_{J-y_0}(y_1)$.  Since the parts $(C_z:z\in V(J)\setminus\{y_0,y_1\})$ are pairwise disjoint and, for each edge $y_1z$ in $J-y_0$, there exists an edge $vw$ of $G$ with $v\in C_{y_1}=N$ and $w\in C_z$, we have $\deg_{J-y_0}(y_1)\le \sum_{v\in N}\deg_G(v) \le \Delta|N| \le \Delta^2 a < n'$, a contradiction. Otherwise, $y_1\in V(Q)$ and $\deg_Q(y_1)\le 3$, so $Q-y_1\subseteq J_2$ contains a $K_{3,n'-1}$-subdivision and hence a $K_{3,n}$-subdivision (as $n'-1\ge n$), again a contradiction. Thus $J-y_0$ contains no $K_{3,n'}$-subdivision.

  Since $\calI$ is a tree-decomposition of $\widehat{I}$ with torsos $\torso{\widehat{I}}{V(I^x)}=I^x$ of treewidth at most $\ell$, $\tw(\widehat{I})\le\ell$. Since $J_2=\widehat{I}-(S\cup N)$ is a subgraph of $\widehat{I}$, $\tw(J_2)\le\ell$. Because $J-y_0$ is obtained by adding the single vertex $y_1$ to $J_2$, $\tw(J-y_0)\le\ell+1$. By \cref{CampbellProduct}, there exist an integer $c\in O((\ell+1)^2(n')^3)=O(\ell^2(n')^3)$ (with $c\ge 1$ depending only on $\ell$ and $n'$, and hence only on $g,p,k,a,\Delta$) and a graph $H$ of treewidth at most $3$ such that $J-y_0$ has an $H$-partition $(D_x:x\in V(H))$ of width at most $c$.  Let $x_1$ be the unique node of $H$ such that $y_1\in D_{x_1}$.  Add a new vertex $x_0$ and the edge $x_0x_1$ to $H$, and set $D_{x_0}:=\{y_0\}$.  Adding the leaf $x_0$ yields $\tw(H)\le\max\{\tw(H-x_0),1\}\le 3$, and since $N_J(y_0)=\{y_1\}$, $(D_x:x\in V(H))$ is an $H$-partition of $J$ of width at most $c$.

  Let $\calP:=(E_x:x\in V(H))$, where $E_x:=\bigcup_{y\in D_x} C_y$ for each $x\in V(H)$. Since $(C_y:y\in V(J))$ is a $J$-partition of $G$ and $(D_x:x\in V(H))$ is an $H$-partition of $J$, $\calP$ is an $H$-partition of $G$. For each $x\in V(H)$ and each $i\in\NN^+$, $|E_x\cap L_i| = \sum_{y\in D_x}|C_y\cap L_i|\le c\ell' \in O(\ell^2(n')^3\ell')$.  Setting $f(g,p,k,a,\Delta):=c\ell'\in\NN^+$, $(\calL,\calP)$ is a $2$-weak layered $H$-partition of $G$ of layered width at most $f(g,p,k,a,\Delta)$ in which $E_{x_0}=C_{y_0}=S\subseteq L_1$ and $\deg_H(x_0)=1$.
\end{proof}

The following lemma is the last big step in the proof of \cref{main}.  The input to this lemma is an almost-embedded tree-decomposition obtained from \cref{dkmw_crisp,promote_apex_neighbours} and the lemma is proven by decomposing the tree $T$ in this tree-decomposition into subtrees that can each be handled by an application of \cref{base_case}.

\begin{lem}\label{full_thing}
  There exists a function $h:\NN^5\to\NN^+$ such that the following is true, for all $g,p,k,a,\Delta\in\NN$:  Let $G$ be a graph of maximum degree at most $\Delta$ and let $(\calT,\calE):=((B_x:x\in V(T)),(\Phi_x:x\in V(T)))$ be a crisp rooted $(g,p,k,a)$-almost-embedded tree-decomposition of $G$ such that for each node $y$ of $T$ with parent $x$, $B_x\cap B_y\subseteq \apex(\Phi_y)$ and $|(B_x\cap B_y)\setminus\topp(\Phi_x)|\le 3$, and let $S\subseteq\apex(\Phi_{r})$, where $r$ is the root of $T$.
  Then there exists a graph $H$ of treewidth at most $3$ and a $2$-weak layered $H$-partition $(\calL:=(L_1,L_2,\ldots),\calP:=(D_x:x\in V(H)))$ of $G$ of layered width at most $h(g,p,k,a,\Delta)$, such that $S=D_{x_0}\subseteq L_1$ for some $x_0\in V(H)$ with $\deg_H(x_0)=1$. 
\end{lem}

\begin{proof}
  For each $x\in V(T)$ and each $u\in\apex(\Phi_x)$, we first bound $|N_{\torso{G}{B_x}}(u)\setminus\topp(\Phi_x)|$.  For each $u\in \apex(\Phi_x)$, $\deg_G(u)\le\Delta$.  For each edge $uv\in E(\torso{G}{B_x})\setminus E(G[B_x])$ with $v\in B_x\setminus\topp(\Phi_x)$, there exists a neighbour $y$ of $x$ in $T$ such that $u,v\in B_x\cap B_y$; since $v\notin\topp(\Phi_x)$ and parent adhesions of $x$ are contained in $\apex(\Phi_x)\subseteq\topp(\Phi_x)$, $y$ must be a child of $x$.  Since $\calT$ is crisp, \cref{crisp_degree} implies that $u$ is in the adhesion $B_x\cap B_y$ for at most $\deg_G(u)\le\Delta$ children $y$ of $x$.  Since $|(B_x\cap B_y)\setminus\topp(\Phi_x)|\le 3$ for each such child $y$ of $x$, this gives at most $3\Delta$ such neighbours $v\in N_{\torso{G}{B_x}}(u)\setminus (N_{G[B_x]}(u)\cup\topp(\Phi_x))$.  Therefore, $|N_{\torso{G}{B_x}}(u)\setminus\topp(\Phi_x)|\le 4\Delta$.

  Let $h(g,p,k,a,\Delta):=f(g+8\Delta a,p+1,k,a,3\Delta)$, where $f$ is the function from \cref{base_case}. Let $(\calT,\calF:=(\Psi_x:x\in V(T)))$ be the rooted $(g+8\Delta a,p+1,k,a)$-almost-embedded tree-decomposition of $G$ obtained by applying \cref{promote_apex_neighbours} with $d:=4\Delta$ to $G$ with the $(g,p,k,a)$-almost-embedded tree-decomposition $(\calT,\calE)$. Thus, for each $x\in V(T)$, $\apex(\Phi_x)=\apex(\Psi_x)$, $\topp(\Phi_x)\subseteq\topp(\Psi_x)$, and $N_{\torso{G}{B_x}}(\apex(\Psi_x))\subseteq \topp(\Psi_x)$.

  The proof is by induction on $|V(T)|$.
  Say that an edge $xy$ of $T$, where $y$ is a child of $x$, is of Type~(i) if $(B_x\cap B_y)\cap\topp(\Psi_x)\neq\emptyset$ and of Type~(ii) if $(B_x\cap B_y)\cap\topp(\Psi_x)=\emptyset$.  
  Refer to \cref{full_thing_fig}.
  Let $T_0$ be the vertex-maximal subtree of $T$ that contains the root of $T$ and whose edges are entirely of Type~(i), and let $T_1,\ldots,T_s$ be the connected components of $T-V(T_0)$. (Possibly $s=0$, in which case $T_0=T$ and we have reached the base of our induction.)  For each $i\in\{1,\ldots,s\}$, let $x_iy_i$ be the unique (Type~(ii)) edge of $T$ with $x_i\in V(T_0)$ and $y_i\in V(T_i)$, root $T_i$ at $y_i$, and let $G_i:=G[\bigcup_{z\in V(T_i)} B_z]$, $\calT_i := (B_x : x \in V(T_i))$, and $\calE_i := (\Phi_x : x \in V(T_i))$.

  \begin{figure}[htbp]
    \begin{center}
      \includegraphics[page=1]{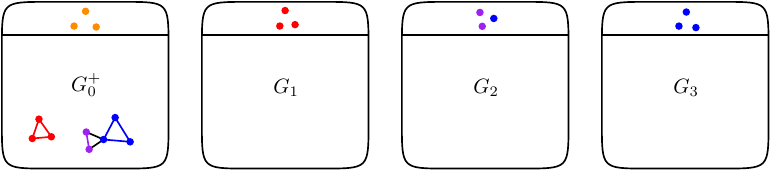} \\      
      \includegraphics[page=2]{figs/full_thing} \\    
      \includegraphics[page=3]{figs/full_thing}     
    \end{center}
    \caption{The main steps in the proof of \cref{full_thing}. The $H_i$-partition of $G_i$ is obtained by induction for each $i\in\{1,\ldots,s\}$. The $H_0$-partition of $G^+_0$ is obtained by applying \cref{base_case}. The final graph $H$ is obtained by making one vertex $x^i_1$ of $H_i-x^i_0$ adjacent to at most 3 vertices of $H_0$ that form a clique, for each $i\in\{1,\ldots,s\}$.}
    \label{full_thing_fig}
  \end{figure}

  For each $i\in\{1,\ldots,s\}$, let $S_i:=B_{x_i}\cap B_{y_i}$.  Since $x_iy_i$ is a Type~(ii) edge, $S_i\cap\topp(\Psi_{x_i})=\emptyset$.  Since $\topp(\Phi_{x_i})\subseteq\topp(\Psi_{x_i})$, we have $S_i\subseteq B_{x_i}\setminus\topp(\Phi_{x_i})$, so $|S_i|=|(B_{x_i}\cap B_{y_i})\setminus\topp(\Phi_{x_i})|\le 3$.  Since $y_i$ is a child of $x_i$ in $T$, $S_i\subseteq\apex(\Phi_{y_i})$.  In particular, $\torso{G_i}{B_{y_i}}-\apex(\Phi_{y_i})=\torso{G}{B_{y_i}}-\apex(\Phi_{y_i})$ and $\torso{G_i}{B_x}=\torso{G}{B_x}$ for all $x\in V(T_i)\setminus\{y_i\}$, so $(\calT_i,\calE_i)$ is a crisp rooted $(g,p,k,a)$-almost-embedded tree-decomposition of $G_i$ satisfying the hypotheses of the lemma (with $S=S_i$).  Since $V(T_i)\subseteq V(T)\setminus V(T_0)\subseteq V(T-r)$, the subtree $T_i$ has fewer vertices than $T$, so we may apply induction on $(\calT_i,\calE_i)$ and $S_i$.
  
  For each $i\in\{1,\ldots,s\}$, apply the inductive hypothesis to $G_i$ with the $(g,p,k,a)$-almost-embedded tree-decomposition $(\calT_i,\calE_i)$ and the vertex set $S_i\subseteq\apex(\Phi_{y_i})$ to obtain a graph $H_i$ of treewidth at most $3$ such that $G_i$ has a $2$-weak layered $H_i$-partition $(\calL_i,\calP_i):=((L^i_1,L^i_2,\ldots),(D_x:x\in V(H_i)))$ of layered width at most $h(g,p,k,a,\Delta)$ in which $S_i=D_{x^i_0}\subseteq L^i_1$ for some $x^i_0\in V(H_i)$ with $\deg_{H_i}(x^i_0)=1$.  For each $i\in\{1,\ldots,s\}$, let $x^i_1$ be the unique neighbour of $x^i_0$ in $H_i$.  For each $i\in\{1,\ldots,s\}$, let $\calW_i:=(W_q:q\in V(Q_i))$ be a tree-decomposition of $H_i-x^i_0$ of width at most $3$. 
  Turn $\calW_i$ into a tree-decomposition of $H_i$ by adding a leaf node $q^i_0$ adjacent to some node $q^i_1$ of $Q_i$ with $x^i_1\in W_{q^i_1}$ and setting $W_{q^i_0}:=\{x^i_0,x^i_1\}$.   

  Let $G_0:=G[\bigcup_{x\in V(T_0)}B_x]$, let $\calT_0:=(B_x:x\in V(T_0))$, and let $\calF_0:=(\Psi_x:x\in V(T_0))$.  Let $G_0^+$ be the graph obtained from $G_0$ by adding edges so that $S_i$ becomes a clique, for each $i\in\{1,\ldots,s\}$. Since $S_i=B_{x_i}\cap B_{y_i}$ is already a clique in $\torso{G}{B_{x_i}}$, $\calT_0$ is a rooted tree-decomposition of $G^+_0$ with $\torso{G^+_0}{B_x}=\torso{G}{B_x}$ for each $x\in V(T_0)$, so $\Psi_x$ is a $(g+8\Delta a,p+1,k,a)$-almost-embedding of $\torso{G^+_0}{B_x}$ satisfying $N_{\torso{G^+_0}{B_x}}(\apex(\Psi_x))\subseteq\topp(\Psi_x)$ for each $x\in V(T_0)$.  Since $\calT$ is crisp, \cref{crisp_degree} implies that each vertex $v\in V(G_0)$ belongs to $S_i=B_{x_i}\cap B_{y_i}$ for at most $\deg_G(v)$ indices $i\in\{1,\ldots,s\}$.  Therefore, for each $v\in V(G_0)$, $\deg_{G_0^+}(v)\le \deg_G(v)+2|\{i\in\{1,\ldots,s\}: v\in S_i\}|\le 3\deg_G(v)\le 3\Delta$.

  By \cref{base_case} applied to the graph $G_0^+$ with the almost-embedded tree-decomposition $(\calT_0,\calF_0)$ and the set $S_0:=S$, there exists a graph $H_0$ of treewidth at most $3$ (chosen so that $V(H_0),V(H_1),\ldots,V(H_s)$ are pairwise disjoint) and a $2$-weak layered $H_0$-partition $(\calL_0,\calP_0):=((L^0_1,L^0_2,\ldots), (D_x:x\in V(H_0)))$ of layered width at most $f(g+8\Delta a,p+1,k,a,3\Delta)=h(g,p,k,a,\Delta)$ in which $S_0=D_{x_0}\subseteq L^0_1$ for some $x_0\in V(H_0)$ with $\deg_{H_0}(x_0)=1$.

  For each $i\in\{1,\ldots,s\}$, there exists an inclusion-minimal set $K_i\subseteq V(H_0)$ of size at most $|S_i|\le 3$ such that $S_i\subseteq\bigcup_{y\in K_i}D_y$.  Since the vertices in $S_i$ form a clique in $G_0^+$, the vertices in $K_i$ form a clique in $H_0$ for each $i\in\{1,\ldots,s\}$.  Since $S\subseteq\apex(\Phi_r)$ and all parent-child adhesions are contained in the apex set of the child, $S\cap B_{x_i}\subseteq\apex(\Phi_{x_i})=\apex(\Psi_{x_i})\subseteq\topp(\Psi_{x_i})$. Because $x_iy_i$ is a Type~(ii) edge of $T$, $S_i\cap\topp(\Psi_{x_i})=\emptyset$, which implies $S\cap S_i=D_{x_0}\cap S_i=\emptyset$, so $x_0\notin K_i$ for each $i\in\{1,\ldots,s\}$. For each $i\in\{1,\ldots,s\}$, merge $H_0$ and $H_i$ by removing $x^i_0$ from $H_i$ and adding the edge $yx^i_1$ for each $y\in K_i$, and call the resulting graph $H$.  Since every edge of $G$ between $S_i$ and $V(G_i)\setminus S_i$ has its endpoints in $\bigcup_{y\in K_i}D_y$ and $D_{x^i_1}$ (as $S_i=D_{x^i_0}$ and $N_{H_i}(x^i_0)=\{x^i_1\}$), $\calP:=(D_y:y\in V(H))$ is an $H$-partition of $G$. Since $x_0\not\in\bigcup_{i=1}^s K_i$, no edges are added incident to $x_0$, so $\deg_H(x_0)=\deg_{H_0}(x_0)=1$.

  To see that the treewidth of $H$ is at most $3$, let $\calW_0:=(W_q:q\in V(Q_0))$ be a tree-decomposition of $H_0$ of width at most $3$.  For each $i\in\{1,\ldots,s\}$, there exists a node $q_i\in V(Q_0)$ such that $K_i\subseteq W_{q_i}$, since the vertices of $K_i$ form a clique in $H_0$. For each $i\in\{1,\ldots,s\}$, join $Q_0$ and $Q_i$ by adding the edge $q_iq^i_0$ and replacing $x^i_0$ with $K_i$ in $W_{q^i_0}$, so that $W_{q^i_0}$ becomes $K_i\cup\{x^i_1\}$.  The resulting tree $Q$ indexes a tree-decomposition $(W_q:q\in V(Q))$ of $H$ of width at most $3$, since $|W_{q^i_0}| = |K_i| + 1 \le 4$.
  
  All that remains is to describe a $2$-weak layering $\calL$ of $G$. 
   For each $i\in\{1,\ldots,s\}$, let $a_i:=1$ if $S_i=\emptyset$, and $a_i:=\min\{j\in\NN^+:L^0_j\cap S_i\neq\emptyset\}$ if $S_i\neq\emptyset$. For each $i\in\{1,\ldots,s\}$, the vertices in $S_i$ form a clique in $G^+_0$, so $S_i\subseteq L^0_{a_i}\cup L^0_{a_i+1}\cup L^0_{a_i+2}$.  For each $v\in V(G_0)$, let $\ell(v)$ be the unique $j\in\NN^+$ such that $v\in L^0_j$.  For each $i\in\{1,\ldots,s\}$ and each $v\in V(G_i)\setminus S_i$, let $\ell(v):=a_i+j-1$ where $j\in\NN^+$ is the unique integer such that $v\in L^i_j$.  Finally, for each $j\in\NN^+$, let $L_j:=\{v\in V(G):\ell(v)=j\}$. 
  
  We now argue that $\calL:=(L_1,L_2,\ldots)$ is a $2$-weak layering of $G$. It is clear that $\calL$ is a cover of $V(G)$ with disjoint sets, so we need only argue that $|\ell(v)-\ell(w)|\le 2$ for each edge $vw$ of $G$. Let $vw$ be an edge of $G$. If $\{v,w\}\subseteq V(G_0)$ then $|\ell(v)-\ell(w)|\le 2$ because $\calL_0$ is a $2$-weak layering of $G^+_0\supseteq G_0$.  If $\{v,w\}\cap V(G_0)=\emptyset$ then there exists a unique $i\in\{1,\ldots,s\}$ such that $\{v,w\}\subseteq V(G_i)\setminus S_i$, so $|\ell(v)-\ell(w)|\le 2$ because $\calL_i$ is a $2$-weak layering of $G_i$.  Otherwise, without loss of generality, $v\in V(G_0)$ and $w\not\in V(G_0)$. Then there exists a unique $i\in\{1,\ldots,s\}$ such that $w\in V(G_i)\setminus S_i$. Thus $v\in S_i=D_{x^i_0}\subseteq L^i_1$. Since $vw\in E(G_i)$ and $\calL_i$ is a $2$-weak layering of $G_i$, $w\in L^i_j$ for some $j\in\{1,2,3\}$, which implies $\ell(w)=a_i+j-1\in\{a_i,a_i+1,a_i+2\}$. Since $v\in S_i\subseteq L^0_{a_i}\cup L^0_{a_i+1}\cup L^0_{a_i+2}$, $\ell(v)\in\{a_i,a_i+1,a_i+2\}$. Since both $\ell(v),\ell(w)\in\{a_i,a_i+1,a_i+2\}$, we have $|\ell(v)-\ell(w)|\le 2$. Hence $\calL$ is a $2$-weak layering of $G$. 
  
  Now, for each $j\in\NN^+$ and each $P\in\calP$, the part $P$ belongs to $\calP_i$ for some $i\in\{0,\ldots,s\}$ and $L_j\cap P$ is contained in at most one layer of $\calL_i$ (namely $L^0_j$ if $i=0$, $L^i_{j-a_i+1}$ if $i\in\{1,\ldots,s\}$ and $j\ge a_i$, or $\emptyset$ if $i\in\{1,\ldots,s\}$ and $j<a_i$).  Therefore, $|L_j\cap P|\le h(g,p,k,a,\Delta)$, so $(\calL,\calP)$ is a $2$-weak layered $H$-partition of $G$ of layered width at most $h(g,p,k,a,\Delta)$, where $H$ is a graph of treewidth at most $3$.  Furthermore, $H$ contains a node $x_0$ of degree $1$ such that $S=D_{x_0}\subseteq L_1$.
\end{proof}

Finally \cref{main} follows easily.

\begin{proof}[Proof of \cref{main}]
  Let $G$ be an $X$-minor-free graph of maximum degree at most $\Delta$. If $V(X)=\emptyset$, then no graph is $X$-minor-free and the statement holds vacuously; assume now that $t:=|V(X)|\ge 1$.  Then $G$ is $K_t$-minor-free. Let $(\calT,\calE)$ be a crisp rooted $(g,p,k,a)$-almost-embedded tree-decomposition of $G$ obtained from \cref{dkmw_crisp}.  Let $r$ be the root of $T$,  and let $S:=\emptyset\subseteq\apex(\Phi_r)$.  By \cref{full_thing}, there exists a $2$-weak layered $H$-partition $(\calL',\calP)$ of $G$ of layered width at most $c:=h(g,p,k,a,\Delta)$, where $H$ is a graph of treewidth at most $3$. By \cref{weak_to_strong}, $G$ has a layered $H$-partition $(\calL,\calP)$ of layered width at most $2c$. Equivalently, $G$ is isomorphic to a subgraph of $H\boxtimes P\boxtimes K_{2c}$ for some path $P$.
\end{proof}

We finish with an open problem. Can the treewidth 3 bound in \cref{main} be improved to treewidth 2? This is open even for planar graphs. Treewidth 1 is impossible, since for any integer $c\geq 1$ there is a planar graph with maximum degree 5 that is not contained in $T\boxtimes P\boxtimes K_c$ where $T$ is a tree and $P$ is a path~\citep{DJMMW24}.  We conjecture that 3 is optimal, even for planar graphs and even for $\Delta=3$:

\begin{conj}\label{main_conjecture}
  For every $c\in\NN^+$ there exists a subcubic planar graph $G$ such that for every graph $H$ of treewidth at most $2$ and every path $P$, the graph $G$ is not isomorphic to any subgraph of $H\boxtimes P\boxtimes K_c$.
\end{conj}

\subsection*{History and Statement of AI Use}

This research began at the MATRIX-hosted workshop \emph{\href{https://www.matrix-inst.org.au/events/global-structure-and-geometry-of-graphs/}{Global Structure and Geometry of Graphs}} held in April 2026 in Creswick, Australia. \Cref{main} was announced and a sketch of the proof was presented at the workshop wrap-up session on April 24, 2026. The current writeup was completed by the authors through emails and shared documents in the subsequent months.

No part of the proof that appears here was generated by AI.  The use of AI was restricted to proofreading, mainly through the use of the Antigravity IDE text editor supported by Gemini.  During the course of this proofreading, Gemini pointed out many typos and ambiguities as well as several gaps in proofs, which were all corrected by hand.

At the time of writing it seems silly to propose a conjecture without first checking if it can be resolved quickly by an LLM.  We have made multiple attempts to ask GPT-6 Astra and Gemini 4 Argon to either disprove \cref{main_conjecture} or prove any weakening of it in which the graph $G$ is $K_t$-minor-free and has maximum-degree $\Delta$ for some fixed $t$ and $\Delta$.  Neither model succeeded, but each produced a long document explaining why they believed it was difficult using current techniques.  In both cases, the models were asked not to return without a resolution in less than five hours.~\textcolor{red}{\faFlag}

{\fontsize{10pt}{11pt}\selectfont
\bibliographystyle{DavidNatbibStyle}
\bibliography{DavidBibliography}}

\end{document}